\documentclass[11pt, a4paper]{article}%
\usepackage{amsmath}
\usepackage{amsfonts}
\usepackage{amssymb}
\usepackage{graphicx}
\usepackage{caption}
\usepackage{subcaption}
\providecommand{\U}[1]{\protect\rule{.1in}{.1in}}
\newtheorem{theorem}{Theorem}[section]

\newtheorem{definition}{Definition}[section]

\newtheorem{lemma}[theorem]{Lemma}
\newtheorem{proposition}[theorem]{Proposition}
\newtheorem{remark}{Remark}[section]

\begin{document}

\title{Optimal ratcheting of dividends in a Brownian risk model}
\author{Hansj\"{o}rg Albrecher\thanks{Department of Actuarial Science, Faculty of
Business and Economics, University of Lausanne, CH-1015 Lausanne and Swiss
Finance Institute. Supported by the Swiss National Science Foundation Project
200021\_191984.}, Pablo Azcue\thanks{Departamento de Matematicas, Universidad
Torcuato Di Tella. Av. Figueroa Alcorta 7350 (C1428BIJ) Ciudad de Buenos
Aires, Argentina.} and Nora Muler$^{\dag}$}
\date{}
\maketitle

\bigskip

\abstract{\begin{quote}
We study the problem of optimal dividend payout from a surplus process governed by Brownian motion with drift under the additional constraint of ratcheting, i.e.\ the dividend rate can never decrease. We solve the resulting two-dimensional optimal control problem, identifying the value function to be the unique viscosity solution of the corresponding Hamilton-Jacobi-Bellman equation. For finitely many admissible dividend rates we prove that threshold strategies are optimal, and for any finite continuum of admissible dividend rates we establish the $\varepsilon$-optimality of curve strategies. This work is a counterpart of \cite{AAM}, where the ratcheting problem was studied for a compound Poisson surplus process with drift. In the present Brownian setup, calculus of variation techniques allow to obtain a much more explicit analysis and description of the optimal dividend strategies. We also give some numerical illustrations of the optimality results.
\end{quote}}

\section{Introduction}

The identification of the optimal way to pay out dividends from a surplus
process to shareholders is a classical topic in actuarial science and
mathematical finance. There is a natural trade-off between paying out gains as
dividends to shareholders early and at the same time leaving sufficient
surplus in order to safeguard future adverse developments and avoid ruin.
Depending on risk preferences, the concrete situation and the simultaneous
exposure to other risk factors such a problem can be formally stated in
various different ways in terms of objective functions and constraints. In
this paper we would like to follow the actuarial tradition of considering the
surplus process as the free capital in an insurance portfolio at any point in
time, and the goal is to maximize the expected sum of discounted dividend
payments that can be paid until the surplus process goes below 0 (which is
called the time of ruin). In such a formulation, the problem goes back to de
Finetti \cite{defin} and Gerber \cite{Ger69}, and has since then been studied
in many variants concerning the nature of the underlying surplus process and
constraints on the type of admissible dividend payment strategies, see
e.g.\ Albrecher \& Thonhauser \cite{AT} and Avanzi \cite{avanzi} for an
overview. From a mathematical perspective, the problem turns out to be quite
challenging, and was cast into the framework of modern stochastic control
theory and the concept of viscosity solutions for corresponding
Hamilton-Jacobi-Bellman equations over the last years, cf.\ Schmidli
\cite{Schmidli book 2008} and Azcue \& Muler \cite{AM Libro}.\newline

Among the variants of the general problem is to look for the optimal dividend
payment strategy if the rate at which dividends are paid can never be reduced.
This \textit{ratcheting} constraint has often been brought up by practitioners
and is in part motivated by the psychological effect that shareholders are
likely to be unhappy about a reduction of dividend payments over time (see
e.g.\ Avanzi et al.\ \cite{aw} for a discussion). One crucial question in this
context is how much of the expected discounted dividends until ruin is lost if
one respects such a ratcheting constraint, if that ratcheting is done in an
optimal way. A first step in that direction was done in Albrecher et
al.\ \cite{ABB}, where the consequences of ratcheting were studied under the
simplifying assumption that a dividend rate can be fixed in the beginning and
can be augmented only once during the lifetime of the process (concretely,
when the surplus process hits some optimally chosen barrier for the first
time). The analysis in that paper was both for a Brownian surplus process as
well as for a surplus process of compound Poisson type. In our recent paper
\cite{AAM}, we then provided the analysis and solution for the general
ratcheting problem for the latter compound Poisson process, and it turned out
that the optimal ratcheting dividend strategy does not lose much efficiency
compared to the unconstrained optimal dividend payout performance, and also
that the one-step ratcheting strategy studied earlier compares remarkably well
to the overall optimal ratcheting solution. In this paper we would like to
address the general ratcheting problem for the Brownian risk model. Such a
model can be seen as a diffusion approximation of the compound Poisson risk
model, but is also interesting in its own right. In particular, the fact that
ruin with zero initial capital is immediate often leads to a more amenable
analysis of stochastic control problems. In addition, the convergence of
optimal strategies from a compound Poisson setting to the one for the
diffusion approximation can be quite delicate, see e.g.\ B\"auerle
\cite{baeuerle}, see also Cohen and Young \cite{cohen} for a recent
convergence rate analysis of simple uncontrolled ruin probabilities towards
its counterparts for the diffusion limit. \newline

On the mathematical level, the general ratcheting formulation leads to a fully
two-dimensional stochastic control problem with all its related challenges,
and it is only recently that in the context of insurance risk theory some
first two-dimensional problems became amenable for analysis, see
e.g.\ Albrecher et al.\ \cite{AAM17}, Gu et al.\ \cite{Gu}, Grandits
\cite{Grandits} and Azcue et al.\ \cite{AMP}. In the present contribution we
would like to exploit the more amenable nature of the ratcheting problem in
the diffusion setting that push the analysis considerably further than was
possible in \cite{AAM}. In particular, we will use calculus of variation
techniques to identify quite explicit formulas for the candidates of optimal
strategies and provide various optimality results. \newline

We would like to mention that optimal ratcheting strategies have been
investigated in the framework of lifetime consumption in the mathematical
finance literature, see e.g.\ Dybvig \cite{dybvig}, Elie and Touzi
\cite{elie}, Jeon et al.\ \cite{Jeon} and more recently Angoshtari et
al.\ \cite{bayr}. However, the concrete model setup and correspondingly also
the involved techniques are quite different to the one of the present paper.
\newline

The remainder of the paper is organized as follows. Section
\ref{Model and basic results} introduces the model and the detailed
formulation of the problem. It also provides some first basic results on
properties of the value function under consideration. Section
\ref{Hamilton-Jacobi-Bellman equations} derives the Hamilton-Jacobi-Bellman
equations and characterization theorems for the value function for both a
closed interval as well as a finite discrete set of admissible dividend
payment rates. In Section \ref{Seccion Convergencia} we prove that the optimal
value function of the problem for discrete sets convergences to the one for a
continuum of admissible dividend rates as the mesh size of the finite set
tends to zero. In Section \ref{Optimal strategies} we show that for finitely
many admissible dividend rates, there exists an optimal strategy for which the
change and non-change regions have only one connected component (this
corresponds to the extension of one-dimensional threshold strategies to the
two-dimensional case). We also provide an implicit equation defining the
optimal threshold function for this case. Subsequently, we turn to the case of
a continuum of admissible dividend rates and use calculus of variation
techniques to identify the optimal curve splitting the state space into a
change and a non-change region as the unique solution of an ordinary
differential equation. We show that the corresponding dividend strategy is
$\varepsilon$-optimal, in the sense that there exists a known sequence of
curves such that the corresponding value functions converge uniformly to the
optimal value function of the problem. Section \ref{Numerical examples}
contains a numerical illustration of the optimal strategy and its performance
relative to the one of for the unconstrained dividend problem and the one
where the dividend rate can only be increased once. Section \ref{seccon}
concludes.\newline

Some technical proofs are delegated to an appendix.

\section{Model and basic results \label{Model and basic results}}

Assume that the surplus process of a company is given by a Brownian motion
with drift%

\[
X_{t}=x+\mu t+\sigma W_{t}%
\]
where $W_{t}$ is a standard Brownian motion, and $\mu,\sigma>0$ are given
constants. Let $(\Omega,\mathcal{F},\left(  \mathcal{F}_{t}\right)  _{t\geq
0},\mathcal{P})$ be the complete probability space generated by the process
$X_{t}$.

The company uses part of the surplus to pay dividends to the shareholders with
rates in a set $S\subset\lbrack0,\overline{c}]$, where $0\leq\overline{c}\in
S$ is the maximum dividend rate possible. Let us denote by $C_{t}$ the rate at
which the company pays dividends at time $t$. Given an initial surplus
$X_{0}=x$ and a minimum dividend rate $c\in S$ at $t=0$, a dividend ratcheting
strategy is given by $C=\left(  C_{t}\right)  _{t\geq0}$ and it is admissible
if it is non-decreasing, right-continuous, adapted with respect to the
filtration $\left(  \mathcal{F}_{t}\right)  _{t\geq0}$ and it satisfies
$C_{t}\in S$ for all $t$. The controlled surplus process can be written as%
\begin{equation}
X_{t}^{C}=X_{t}-\int_{0}^{t}C_{s}ds. \label{XtC}%
\end{equation}
Define $\Pi_{x,c}^{S}$ as the set of all admissible dividend ratcheting
strategies with initial surplus $x\geq0$ and minimum initial dividend rate
$c\in S$. \ Given $C\in\Pi_{x,c}^{S}$, the value function of this strategy is
given by%
\[
J(x;C)=\mathbb{E}\left[  \int_{0}^{\tau}e^{-qs}C_{s}ds\right]  ,
\]
where $q>0$ and $\tau=\inf\left\{  t\geq0:X_{t}^{C}<0\right\}  $ is the ruin
time. Hence, for any initial surplus $x\geq0$ and initial dividend rate $c$,
our aim is to maximize%

\begin{equation}
V^{S}(x,c)=\sup_{C\in\Pi_{x,c}^{S}}J(x;C). \label{Optimal Value Function}%
\end{equation}

It is immediate to see that $V^{S}(0,c)=0$ for all $c\in S.$

\begin{remark}
\normalfont
\label{Optima sin ratcheting} The dividend optimization problem without the
ratcheting constraint, that is where the dividend strategy $C=\left(
C_{t}\right)  _{t\geq0}$ is not necessarily non-decreasing, was studied
intensively in the literature (see e.g.\ Gerber and Shiu \cite{Gerber2004}).
Unlike the ratcheting optimization problem, this non-ratcheting problem is
one-dimensional. If $V_{NR}(x)$ denotes its optimal value function, then
clearly $V^{S}(x,c)\leq V_{NR}(x)$ for all $x\geq0$ and $c\in S\subset
\lbrack0,\overline{c}]$. The function $V_{NR}$ is increasing, concave, twice
continuously differentiable with $V_{NR}(0)=0$ and $\lim_{x\rightarrow\infty
}V_{NR}(x)=\overline{c}/q$; so it is Lipschitz with Lipschitz constant
$V_{NR}^{\prime}(0).$
\end{remark}

The following Lemma states the dynamic programming principle, its proof is
similar to the one of Lemma 1.2 in Azcue and Muler \cite{AM Libro}{\Large .}

\begin{lemma}
\label{DPP} Given any stopping time $\widetilde{\tau}$, we can write
\[
V^{S}(x,c)=\sup\limits_{C\in\Pi_{x,c}^{S}}\mathbb{E}\left[  \int_{0}
^{\tau\wedge\widetilde{\tau}}e^{-qs}C_{s}ds+e^{-q(\tau\wedge\widetilde{\tau}
)}V^{S}(X_{\tau\wedge\widetilde{\tau}}^{C},C_{\tau\wedge\widetilde{\tau}
})\right]  \text{.}
\]

\end{lemma}

\bigskip

We now state a straightforward result regarding the boundedness and
monotonicity of the optimal value function.

\begin{proposition}
\label{Monotone Optimal Value Function}The optimal value function $V^{S}(x,c)
$ is bounded by $\overline{c}/q$, non-decreasing in $x$ and non-increasing in
$c.$
\end{proposition}

\textit{Proof.} Since the discounted value of paying the maximum rate
$\overline{c}$ up to infinity is $\overline{c}/q,$ we conclude the boundedness result.

On the one hand $V^{S}(x,c)$ is non-increasing in $c$ because given
$c_{1}<c_{2}$ we have $\Pi_{x,c_{2}}^{S}\subset\Pi_{x,c_{1}}^{S}$ for any
$x\geq0$. On the other hand, given $x_{1}<x_{2}$ and an admissible ratcheting
strategy $C^{1}\in\Pi_{x_{1},c}^{S}$ for any $c\in S$, let us define $C^{2}%
\in\Pi_{x_{2},c}^{S}$ as $C_{t}^{2}=C_{t}^{1}$ until the ruin time of the
controlled process $X_{t}^{C^{1}}$ with $X_{0}^{C^{1}}=x_{1}$, and pay the
maximum rate $\overline{c}$ afterwards. Thus, $J(x;C_{1})\leq J(x;C_{2})$ and
we have the result.\hfill$\blacksquare$\newline

The Lipschitz property of the function $V_{NR}$ introduced in Remark
\ref{Optima sin ratcheting} can now be used to prove a first result on the
regularity of the function $V^{S}.$

\begin{proposition}
\label{Proposition Global Lipschitz zone}There exists a constant $K>0$ such
that
\[
0\leq V^{S}(x_{2},c_{1})-V^{S}(x_{1},c_{2})\leq K\left[  \left(  x_{2}
-x_{1}\right)  +\left(  c_{2}-c_{1}\right)  \right]
\]

for all $0\leq x_{1}\leq x_{2}$ and $c_{1},c_{2}\in S$ with $c_{1}\leq c_{2}.
$
\end{proposition}

The proof is given in the appendix.

\section{Hamilton-Jacobi-Bellman equations
\label{Hamilton-Jacobi-Bellman equations}}

In this section we introduce the Hamilton-Jacobi-Bellman (HJB) equation of the
ratcheting problem for $S\subset\lbrack0,\infty)$, when $S$ is either a closed
interval or a finite set. We show that the optimal value function $V$ defined
in (\ref{Optimal Value Function}) is the unique viscosity solution of the
corresponding HJB equation with boundary condition $\overline{c}/q$ when $x$
goes to infinity, where $\overline{c}=\max S.$

First, consider the case $S=\left\{  c\right\}  $. In this case, the unique
admissible strategy consists of paying a constant dividend rate $c $ up to the
ruin time. Correspondingly, the value function $V^{\left\{  c\right\}  }(x,c)$
is the unique solution of the second order differential equation%
\begin{equation}
\mathcal{L}^{c}(W):=\frac{\sigma^{2}}{2}\partial_{xx}W+(\mu-c)\partial
_{x}W-qW+c=0 \label{Lc}%
\end{equation}
with boundary conditions $V^{\left\{  c\right\}  }(0,c)=0$ and $\lim
_{x\rightarrow\infty}$ $V^{\left\{  c\right\}  }(x,c)=c/q.$ The solutions of
\eqref{Lc} are of the form
\begin{equation}
\frac{c}{q}+a_{1}e^{\theta_{1}(c)x}+a_{2}e^{\theta_{2}(c)x}\text{ with }
a_{1},a_{2}\in{\mathbb{R}} , \label{Solucion General L=0}%
\end{equation}
where $\theta_{1}(c)>0$ and $\theta_{2}(c)<0$ are the roots of the
characteristic equation
\[
\frac{\sigma^{2}}{2}z^{2}+(\mu-c)z-q=0
\]
associated to the operator $\mathcal{L}^{c}$, that is%
\begin{equation}
\theta_{1}(c):=\frac{c-\mu+\sqrt{(c-\mu)^{2}+2q\sigma^{2}}}{\sigma^{2}}%
,\quad\theta_{2}(c):=\text{$\frac{c-\mu-\sqrt{(c-\mu)^{2}+2q\sigma^{2}}}
{\sigma^{2}}$.} \label{Definicion tita1 tita2}%
\end{equation}

In the following remark, we state some basic properties of $\theta_{1}$ and
$\theta_{2}.$

\begin{remark}
\normalfont
\label{Propiedades de Titas} We have that

\begin{enumerate}
\item $\theta_{1}(c)=-\theta_{2}(c)$ if $c=\mu$ and $\theta_{1}^{2}
(c)\geq\theta_{2}^{2}$$(c)$ if, and only if, $c-\mu\geq0.$

\item $\theta_{1}^{\prime}(c)=\frac{1}{\sigma^{2}}(1+\frac{c-\mu}{\sqrt
{(c-\mu)^{2}+2q\sigma^{2}}})$ and $\theta_{2}^{\prime}(c)=\frac{1}{\sigma^{2}
}(1-\frac{c-\mu}{\sqrt{(c-\mu)^{2}+2q\sigma^{2}}}),$ so $\theta_{1}^{\prime
}(c),\theta_{2}^{\prime}(c)\in(0,\frac{2}{\sigma^{2}}).$
\end{enumerate}
\end{remark}

\noindent The solutions of $\mathcal{L}^{c}(W)=0$ with boundary condition
$W(0)=0$ are of the form%

\begin{equation}
\frac{c}{q}\left(  1-e^{\theta_{2}(c)x}\right)  +a(e^{\theta_{1}
(c)x}-e^{\theta_{2}(c)x})\ \text{with }a\in{\mathbb{R}}.
\label{Solucion de L=0 con condicion en 0}%
\end{equation}
And finally, the unique solution of $\mathcal{L}^{c}(W)=0$ with boundary
conditions $W(0)=0$ and $\lim_{x\rightarrow\infty}$ $W(x)=c/q\ $corresponds to
$a=0,$ so that%
\begin{equation}
V^{\left\{  c\right\}  }(x,c)=\frac{c}{q}\left(  1-e^{\theta_{2}(c)x}\right)
. \label{Formula de V con c constante}%
\end{equation}
We have that $V^{\left\{  c\right\}  }(\cdot,c)$ is increasing and concave.

\begin{remark}
\normalfont
\label{lim l0/c}Given a set $S\subset\lbrack0,\infty)$ with $\overline{c}=\max
S<\infty,$ we have that
\[
V^{S}(x,c)\geq V^{\left\{  \overline{c}\right\}  }(x,\overline{c}
)=\frac{\overline{c}}{q}\left(  1-e^{\theta_{2}(\overline{c})x}\right)
\]
and so, by Remark \ref{Optima sin ratcheting}, we conclude that $\lim
_{x\rightarrow\infty}V^{S}(x,c)={\overline{c}}/{q}$ for any $c\in S$.
\end{remark}

\subsection{Hamilton-Jacobi-Bellman equations for closed intervals
\label{Hamilton-Jacobi-Bellman equations for closed intervals}}

Let us now consider the case $S=[\underline{c},\overline{c}]$ with
$0\leq\underline{c}<\overline{c}.$ The HJB equation associated to
(\ref{Optimal Value Function}) is given by%

\begin{equation}
\max\{\mathcal{L}^{c}(u)(x,c),\partial_{c}u(x,c)\}=0\text{ for }
x\geq0\ \text{and }\underline{c}\leq c\leq\overline{c}\text{.}
\label{HJB equation}%
\end{equation}

We say that a function $f:[0,\infty)\times\lbrack\underline{c},\overline
{c})\rightarrow{\mathbb{R}}$ is \textit{(2,1)-differentiable} if $f$ is
continuously differentiable and $\partial_{x}f(\cdot,c)$ is continuously differentiable.

\bigskip

\begin{definition}
\label{Viscosity}

(a) A locally Lipschitz function $\overline{u}:[0,\infty)\times\lbrack
\underline{c},\overline{c})\rightarrow{\mathbb{R}}$\ where $0\leq\underline
{c}<\overline{c}$ is a viscosity supersolution of (\ref{HJB equation})\ at
$(x,c)\in(0,\infty)\times\lbrack\underline{c},\overline{c})$,\ if any
(2,1)-differentiable function $\varphi:[0,\infty)\times\lbrack\underline
{c},\overline{c})\rightarrow{\mathbb{R}}\ $with $\varphi(x,c)=\overline
{u}(x,c)$ such that $\overline{u}-\varphi$\ reaches the minimum at $\left(
x,c\right)  $\ satisfies
\[
\max\left\{  \mathcal{L}^{c}(\varphi)(x,c),\partial_{c}\varphi(x,y)\right\}
\leq0.\
\]
The function $\varphi$ is called a \textbf{test function for supersolution} at
$(x,c)$.

(b) A function $\underline{u}:$ $[0,\infty)\times\lbrack\underline
{c},\overline{c})\rightarrow{\mathbb{R}}\ $\ is a viscosity subsolution\ of
(\ref{HJB equation})\ at $(x,c)\in(0,\infty)\times\lbrack\underline
{c},\overline{c})$,\ if any (2,1)-differentiable function $\psi:[0,\infty
)\times\lbrack\underline{c},\overline{c})\rightarrow{\mathbb{R}}\ $with
$\psi(x,c)=\underline{u}(x,c)$ such that $\underline{u}-\psi$\ reaches the
maximum at $\left(  x,c\right)  $ satisfies
\[
\max\left\{  \mathcal{L}^{c}(\psi)(x,c),\partial_{c}\psi(x,c)\right\}
\geq0\text{.}
\]
The function $\psi$ is called a \textbf{test function for subsolution} at
$(x,c)$.

(c) A function $u:[0,\infty)\times\lbrack\underline{c},\overline
{c})\rightarrow{\mathbb{R}}$ which is both a supersolution and subsolution at
$(x,c)\in\lbrack0,\infty)\times\lbrack\underline{c},\overline{c})$ is called a
viscosity solution of (\ref{HJB equation})\ at $(x,c)$.
\end{definition}

\begin{remark}
\normalfont
\label{No importa el c}Note that, by (\ref{Optimal Value Function}),
$V^{[0,\overline{c}]}(x,c)=V^{[\underline{c},\overline{c}]}(x,c)$ for all
$0\leq\underline{c}\leq c\leq\overline{c}$, so in order to simplify the
notation we define $V(x,c):=V^{[\underline{c},\overline{c}]}(x,c):[0,\infty
)\times\lbrack\underline{c},\overline{c})\rightarrow{\mathbb{R}.}$
\end{remark}

We first prove that $V$ is a viscosity solution of the corresponding HJB equation.

\begin{proposition}
\label{Proposicion Viscosidad} $V$ is a viscosity solution of
(\ref{HJB equation}) in $(0,\infty)\times\lbrack\underline{c},\overline{c})$.
\end{proposition}

The proof is given in the appendix. \newline

Note that by definition of ratcheting $V(x,\overline{c})$ corresponds to the
value function of the strategy that constantly pays dividends at rate
$\overline{c}$, with initial surplus $x$. So, by
(\ref{Formula de V con c constante}),
\begin{equation}
V(x,\overline{c})=V^{\{\overline{c}\}}(x,\overline{c}). \label{V en cbarra}%
\end{equation}

Let us now state the comparison result for viscosity solutions.

\begin{lemma}
\label{Lema para Unicidad} Assume that (i) $\underline{u}$ is a viscosity
subsolution and $\overline{u}$ is a viscosity supersolution of the HJB
equation (\ref{HJB equation}) for all $x>0$ and for all $c\in\lbrack
\underline{c},\overline{c})$ with $0\leq\underline{c}<\overline{c}$, (ii)
$\underline{u}$ and $\overline{u}$ are non-decreasing in the variable $x$ and
Lipschitz
in $[0,\infty)\times\lbrack\underline{c},\overline{c})$, and (iii)
$\underline{u}(0,c)=\overline{u}(0,c)=0$, $\lim_{x\rightarrow\infty}%
\underline{u}(x,c)\leq\overline{c}/q\leq\lim_{x\rightarrow\infty}\overline
{u}(x,c)$. Then $\underline{u} \leq\overline{u}$ in $[0,\infty)\times
\lbrack\underline{c},\overline{c}).$
\end{lemma}

The proof is given in the appendix. \newline

The following characterization theorem is a direct consequence of the previous
lemma, Remark \ref{lim l0/c} and Proposition \ref{Proposicion Viscosidad}.

\begin{theorem}
\label{Caracterizacion Continua}The optimal value function $V$ is the unique
function non-decreasing in $x$ that is a viscosity solution of
(\ref{HJB equation}) in $(0,\infty)\times\lbrack\underline{c},\overline{c})$
with $V(0,c)=0$ and $\lim_{x\rightarrow\infty}$ $V(x,c)=\overline{c}/q$ for
$c\in\lbrack\underline{c},\overline{c}).$
\end{theorem}

From Definition \ref{Optimal Value Function}, Lemma \ref{Lema para Unicidad},
and Remark \ref{lim l0/c} together with Proposition
\ref{Proposicion Viscosidad}, we also get the following verification theorem.

\begin{theorem}
\label{verification result} Consider $S=[\underline{c},\overline{c}]$ and
consider a family of strategies
\[
\left\{  C_{x,c}\in\Pi_{x,c}^{S}:(x,c)\in\lbrack0,\infty)\times\lbrack
\underline{c},\overline{c}]\right\}  .
\]
If the function $W(x,c):=J(x;C_{x,c})$ is a viscosity supersolution of the HJB
equation (\ref{HJB equation}) in $(0,\infty)\times\lbrack\underline
{c},\overline{c})$ with $W(0,c)=0$ and $\lim_{x\rightarrow\infty}W(x,c)=$
$\overline{c}/q,$ then $W$ is the optimal value function $V$. Also, if for
each $k\geq1$ there exists a family of strategies $\left\{  C_{x,c}^{k}\in
\Pi_{x,c}^{S}:(x,c)\in\lbrack0,\infty)\times\lbrack\underline{c},\overline
{c}]\right\}  $ such that $W(x,c):=\lim_{k\rightarrow\infty}J(x;C_{x,c}^{k})$
is a viscosity supersolution of the HJB equation (\ref{HJB equation}) in
$(0,\infty)\times\lbrack\underline{c},\overline{c})$ with $W(0,c)=0$ and
$\lim_{x\rightarrow\infty}W(x,c)=$ $\overline{c}/q$, then $W$ is the optimal
value function $V$.
\end{theorem}

\subsection{Hamilton-Jacobi-Bellman equations for finite sets
\label{Hamilton-Jacobi-Bellman equations for finite sets}}

Let us now consider the case%
\[
S=\left\{  c_{1},c_{2},....,c_{n}\right\}  ,
\]
where $0\leq c_{1}<c_{2}<....<c_{n}=\overline{c}$. Note that $V^{S}%
(x,c_{i})=V^{\left\{  c_{i},c_{i+1},....,c_{n}\right\}  }(x,c_{i}).$ We
simplify the notation as follows:%

\begin{equation}
V^{c_{i}}(x):=V^{S}(x,c_{i}). \label{optimal value function Discreta}%
\end{equation}
So we have the inequalities
\[
V^{c_{i}}(x)\geq V^{c_{i+1}}(x)\geq...\geq V^{c_{n}}(x)=V^{\overline{c}%
}(x)\geq0,
\]
where $V^{\overline{c}}(x)=V^{\{\overline{c}\}}(x,\overline{c})$ as defined in
(\ref{Formula de V con c constante}).

The Hamilton-Jacobi-Bellman equation associated to
(\ref{optimal value function Discreta}) is given by%
\begin{equation}
\max\left\{  \mathcal{L}^{c_{i}}(V^{c_{i}}(x)),V^{c_{i+1}}(x)-V^{c_{i}
}(x)\right\}  =0\text{ for }x\geq0\text{ and }i=1,...,n-1\text{.}
\label{HJB equation discreta}%
\end{equation}

As in the continuous case we have that $V^{c_{i}}$ is the viscosity solution
of the corresponding HJB equation. Let us introduce the definition of a
viscosity solution in the one-dimensional case.

\begin{definition}
\bigskip(a) A locally Lipschitz function $\overline{u}:[0,\infty
)\rightarrow{\mathbb{R}}$\ is a viscosity supersolution of
(\ref{HJB equation discreta})\ at $x\in(0,\infty)$\ if any twice continuously
differentiable function $\varphi:[0,\infty)\rightarrow{\mathbb{R}}\ $with
$\varphi(x)=\overline{u}(x)$ such that $\overline{u}-\varphi$\ reaches the
minimum at $x$\ satisfies
\[
\max\left\{  \mathcal{L}^{c_{i}}(\varphi)(x),V^{c_{i+1}}(x)-\varphi
(x)\right\}  \leq0.\
\]
The function $\varphi$ is called a \textbf{test function for supersolution} at
$x$.

(b) A function $\underline{u}:$ $[0,\infty)\rightarrow{\mathbb{R}}\ $\ is a
viscosity subsolution\ of (\ref{HJB equation discreta})\ at $x\in(0,\infty
)$\ if any twice continuously differentiable function $\psi:[0,\infty
)\rightarrow{\mathbb{R}}\ $with $\psi(x)=\underline{u}(x)$ such that
$\underline{u}-\psi$\ reaches the maximum at $x$ satisfies
\[
\max\left\{  \mathcal{L}^{c_{i}}(\psi)(x),V^{c_{i+1}}(x)-\psi(x)\right\}
\geq0\text{.}
\]
The function $\psi$ is called a \textbf{test function for subsolution} at $x$.

(c) A function $u:[0,\infty)\rightarrow{\mathbb{R}}$ which is both a
supersolution and subsolution at $x\in\lbrack0,\infty)$ is called a viscosity
solution of (\ref{HJB equation discreta})\ at $x$.
\end{definition}

The following characterization theorem is the analogue of Theorem
\ref{Caracterizacion Continua} for finite sets; the proof is similar and
simpler than the one in the continuous case.

\begin{theorem}
\label{Caracterizacion Discreta}The optimal value function $V^{c_{i}}(x)$ for
$1\leq i<n$ is the unique viscosity solution of the associated HJB equation
(\ref{HJB equation discreta}) with boundary condition $V^{c_{i}}(0)=0$ and
$\lim_{x\rightarrow\infty}V^{c_{i}}(x)=\overline{c}/q.$
\end{theorem}

We also have the alternative characterization theorem.

\begin{theorem}
\label{Menor supersolucion discreta}The optimal value function $V^{c_{i}}(x)$
for $1\leq i<n$ is the smallest viscosity supersolution of the the associated
HJB equation (\ref{HJB equation discreta}) with boundary condition $0$ at
$x=0$ and limit greater than or equal to $\overline{c}/q\ $as $x$ goes to infinity.
\end{theorem}

\begin{remark}
\normalfont
\label{remark del obstaculo}The function $V^{c_{n}}$ has the closed formula
given by (\ref{Formula de V con c constante}) for $c=c_{n}.$ By the previous
theorem, once $V^{c_{i+1}}$ is known, the optimal value function $V^{c_{i}}$
can be obtained recursively as the solution of the \textbf{obstacle problem}
of finding the smallest viscosity supersolution of the equation $\mathcal{L}
^{c_{i}}=0$ above the obstacle $V^{c_{i+1}}$.
\end{remark}

\section{Convergence of the optimal value functions from the discrete to the
continuous case \label{Seccion Convergencia}}

In this section we prove that the optimal value functions of finite ratcheting
strategies approximate the optimal value function of the continuous case as
the mesh size of the finite sets goes to zero.

Consider for $n\geq0$, a sequence of sets $\mathcal{S}^{n}$ (with $k_{n}$
elements) of the form{\small \ }
\[
\mathcal{S}^{n}=\left\{  c_{1}^{n}=\underline{c}<c_{2}^{n}<\cdots<c_{k_{n}%
}^{n}=\overline{c}\right\}  .
\]
satisfying $\mathcal{S}^{0}=\left\{  \underline{c},\overline{c}\right\}  $,
$\mathcal{S}^{n}\subset\mathcal{S}^{n+1}$ and mesh-size $\delta(\mathcal{S}%
^{n}):=\max_{k=2,k_{n}}\left(  c_{k}^{n}-c_{k-1}^{n}\right)  \searrow0$ as $n$
goes to infinity.

Let us extend the definition of $V^{\mathcal{S}^{n}}$ to the function
$V^{n}:[\underline{c},\infty)\times\lbrack0,\overline{c}]\rightarrow
{\mathbb{R}},$ as%
\begin{equation}
V^{n}(x,c)=V^{\mathcal{S}^{n}}(x,\widetilde{c}^{n}), \label{Definicion Vn}%
\end{equation}
where
\begin{equation}
\widetilde{c}^{n}=\min\{c_{i}^{n}\in\mathcal{S}^{n}:c_{i}^{n}\geq c\}.
\label{Deminicioncruliton}%
\end{equation}
We will prove that $\lim_{n\rightarrow\infty}V^{n}(x,c)=V^{[\underline
{c},\overline{c}]}(x,c)$ for any $(x,c)\in\lbrack0,\infty)\times
\lbrack\underline{c},\overline{c}]$ and we will study the uniform convergence
of this limit.

Since $V^{n}\leq V^{n+1}\leq V^{[\underline{c},\overline{c}]}$, there exists
the limit function
\begin{equation}
\overline{V}(x,c):=\lim_{n\rightarrow\infty}V^{n}(x,c).
\label{Definicion Vbarra}%
\end{equation}
Later on, we will show that $\overline{V}=V^{[\underline{c},\overline{c}]}$.
Note that $\overline{V}(x,c)$ is non-increasing in $c$ with $\overline
{V}(x,\overline{c})=V(x,\overline{c})$, and non-decreasing in $x$ with
$\lim_{x\rightarrow\infty}$ $\overline{V}(x,c)=\overline{c}/q$. With the same
proof the one for Proposition 6.1 of \cite{AAM}, we have the following proposition:

\begin{proposition}
\label{Proposicion Convergencia Uniforme VN a V}The sequence $V^{n}$ converges
uniformly to $\overline{V}.$
\end{proposition}

With this, we can obtain the main result of this section.

\begin{theorem}
\label{Vbarra es V} The function $\overline{V}$ defined in
(\ref{Definicion Vbarra}) is the optimal value function $V^{[\underline
{c},\overline{c}]}$.
\end{theorem}

\textit{Proof}. Note that $\overline{V}(x,c)$ is a limit of value functions of
admissible strategies, so in order to satisfy the assumptions of Theorem
\ref{verification result}, it remains to see that $\overline{V}$ is a
viscosity supersolution of (\ref{HJB equation}) at any point $(x_{0},c_{0})$
with $x_{0}>0.$ $\partial_{c}\overline{V}(x_{0},c_{0})\leq0$ in the viscosity
sense because $\overline{V}$ is non-increasing in $c$; so it is sufficient to
show that $\mathcal{L}^{c_{0}}(\overline{V})(x_{0},c_{0})\leq0$ in the
viscosity sense. Let $\varphi$ be a test function for viscosity supersolution
of (\ref{HJB equation}) at $(x_{0},c_{0})$, i.e.\ a (2,1)-differentiable
function $\varphi$ with
\begin{equation}
\overline{V}(x,c)\geq\varphi(x,c)\text{ and }\overline{V}(x_{0},c_{0}
)=\varphi(x_{0},c_{0})\text{.} \label{Comparacion}%
\end{equation}
In order to prove that $\mathcal{L}^{c}(\varphi)(x_{0},c_{0})\leq0$, consider
now, for $\gamma>0$ small enough,
\[
\varphi_{\gamma}(x,c)=\varphi(x,c)-\gamma(x-x_{0})^{4}.
\]
Given $n\geq1,$ let us consider $\widetilde{c}_{0}^{n}$ as defined in
(\ref{Deminicioncruliton}),%
\[
a_{n}^{\gamma}:=\min\{V^{n}(x,\widetilde{c}_{0}^{n})-\varphi_{\gamma
}(x,\widetilde{c}_{0}^{n}):x\in\lbrack0,x_{0}+1]\},
\]%
\[
x_{n}^{\gamma}:=\arg\min\{V^{n}(x,\widetilde{c}_{0}^{n})-\varphi_{\gamma
}(x,\widetilde{c}_{0}^{n}):x\in\lbrack0,x_{0}+1]\},
\]
and%
\[
b_{n}^{\gamma}:=\max\{\overline{V}(x,\widetilde{c}_{0}^{n})-V^{n}%
(x,\widetilde{c}_{0}^{n}):x\in\lbrack0,x_{0}+1]\}.
\]
Since $\widetilde{c}_{0}^{n}\searrow c_{0}$ and, from Proposition
\ref{Proposicion Convergencia Uniforme VN a V}, $\lim_{n\rightarrow\infty
}a_{n}^{\gamma}=0$ and $\lim_{n\rightarrow\infty}b_{n}^{\gamma}=0$, we also
have that $\lim_{n\rightarrow\infty}~x_{n}^{\gamma}=x_{0}$ because%
\[%
\begin{array}
[c]{lll}%
0 & = & V^{n}(x_{n}^{\gamma},\widetilde{c}_{0}^{n})-\left(  \varphi_{\gamma
}(x_{n}^{\gamma},\widetilde{c}_{0}^{n})+a_{n}^{\gamma}\right) \\
& = & \left(  V^{n}(x_{n}^{\gamma},\widetilde{c}_{0}^{n})-\overline{V}
(x_{n}^{\gamma},\widetilde{c}_{0}^{n})\right)  +\left(  \overline{V}
(x_{n}^{\gamma},\widetilde{c}_{0}^{n})-\varphi_{\gamma}(x_{n}^{\gamma
},\widetilde{c}_{0}^{n})\right)  -a_{n}^{\gamma}\\
& \geq & -b_{n}^{\gamma}+0-a_{n}^{\gamma}+\gamma(x_{n}^{\gamma}-x_{0})^{4}%
\end{array}
\]
and then%
\[
(x_{n}^{\gamma}-x_{0})^{4}\leq\frac{a_{n}^{\gamma}+b_{n}^{\gamma}}{\gamma
}\rightarrow0~\text{as }n\rightarrow\infty.
\]

Note that $\overline{\varphi}^{n}(\cdot)=\varphi_{\gamma}(\cdot,\widetilde
{c}_{0}^{n})+a_{n}^{\gamma}$ is a test function for viscosity supersolution of
$V^{n}(\cdot,\widetilde{c}_{0}^{n})$ in Equation (\ref{HJB equation discreta})
at the point $x_{n}^{\gamma}$ because%

\[
\varphi_{\gamma}(x_{n}^{\gamma},\widetilde{c}_{0}^{n})+a_{n}^{\gamma}%
=V^{n}(x_{n}^{\gamma},\widetilde{c}_{0}^{n})\text{ and }\varphi_{\gamma
}(x,\widetilde{c}_{0}^{n})+a_{n}^{\gamma}\leq V^{n}(x,\widetilde{c}_{0}%
^{n})\text{ for }x\in\lbrack0,x_{0}+1].
\]
And so%
\[
\mathcal{L}^{\widetilde{c}_{0}^{n}}(\varphi_{\gamma})(x_{n}^{\gamma
},\widetilde{c}_{0}^{n})=\mathcal{L}^{\widetilde{c}_{0}^{n}}(\overline
{\varphi}^{n})(x_{n}^{\gamma})+qa_{n}^{\gamma}\leq qa_{n}^{\gamma}.
\]
Since $(x_{n}^{\gamma},c_{n})\rightarrow(x_{0},c_{0})$, $\overline{\varphi
}^{n}(\cdot)=\varphi_{\gamma}(\cdot,\widetilde{c}_{0}^{n})+a_{n}^{\gamma
}\rightarrow\varphi_{\gamma}(\cdot,c_{0})$ as $n\rightarrow\infty$ and
$\varphi_{\gamma}$ is (2,1)-differentiable, one gets
\[
\mathcal{L}^{c_{0}}(\varphi_{\gamma})(x_{0},c_{0})=\lim_{n\rightarrow\infty
}\mathcal{L}^{\widetilde{c}_{0}^{n}}(\overline{\varphi}^{n})(x_{n}^{\gamma
})\leq0.
\]
Finally, as
\[
\partial_{x}\varphi_{\gamma}(x_{0},c_{0})=\partial_{x}\varphi(x_{0}%
,c_{0})~\text{and }\partial_{xx}\varphi_{\gamma}(x_{0},c_{0})=\partial
_{xx}\varphi(x_{0},c_{0})
\]
and $\varphi_{\gamma}\nearrow\varphi$ as $\gamma\searrow0$, we obtain that
$\mathcal{L}^{c_{0}}(\varphi)(x_{0},c_{0})\leq0$ and the result follows.\hfill
$\blacksquare$

\section{The optimal strategies \label{Optimal strategies}}

We show first that, regardless whether $S$ is finite or an interval with $\max
S=\overline{c}$, the optimal strategy for sufficiently small $\overline{c}$ is
to immediately start paying dividends at the maximum rate $\overline{c}$.

\begin{proposition}
\label{Estrategia trivial} If $\overline{c}\leq q\sigma^{2}/(2\mu)$, then
$V(x,c)=V^{\{\overline{c}\}}(x,\overline{c})$ for any $(x,c)\in\lbrack
0,\infty)\times S.$
\end{proposition}

\textit{Proof.} If we call $W(x,c):=V^{\{\overline{c}\}}(x,\overline{c}),$
then we know that $\mathcal{L}^{\overline{c}}(W)(x,c)=0$. Since $W(0,\overline
{c})=0,$ $\lim_{x\rightarrow\infty}W(x,c)=\overline{c}/q$ and $\partial
_{c}W(x,c)=0$, then by Theorem \ref{Caracterizacion Continua} and Theorem
\ref{Caracterizacion Discreta} it is enough to prove that $\mathcal{L}%
^{c}(W)(x,c)\leq0$ for $c\in S.$ But, by (\ref{Formula de V con c constante})
and (\ref{Definicion tita1 tita2})
\[%
\begin{array}
[c]{lll}%
\mathcal{L}^{c}(W)(x,c) & = & \mathcal{L}^{\overline{c}}(W)(x,c)+(\overline
{c}-c)(\partial_{x}W(x,c)-1)\\
& = & (\overline{c}-c)(-\frac{\overline{c}}{q}\theta_{2}(\overline
{c})e^{\theta_{2}(\overline{c})x}-1)\\
& \leq & (\overline{c}-c)(-\frac{\overline{c}}{q}\theta_{2}(\overline{c})-1)\\
& \leq & 0
\end{array}
\]
for $c\leq\overline{c}\leq\frac{q\sigma^{2}}{2\mu}$.\hfill$~\blacksquare$

\begin{remark}
\normalfont
\label{Trivial en Entorno} The proof of the previous proposition also shows
that if there exists a $d\in$ $S\setminus\{\overline{c}\}$ and $\varepsilon>0$
such that $V(x,d)=V^{\{\overline{c}\}}(x,\overline{c})$ for $x\in
\lbrack0,\varepsilon]$, then $\overline{c}\leq\frac{q\sigma^{2}}{2\mu}$ and so
$V(x,c)=V^{\{\overline{c}\}}(x,\overline{c})$ for any $(x,c)\in\lbrack
0,\infty)\times S.$
\end{remark}

Let us first address the case of $S=\left\{  c_{1},c_{2},....,c_{n}\right\}  $
with $0\leq c_{1}<c_{2}<....<c_{n}=\overline{c}$. We introduce the concept of
strategies with a threshold structure for each level $c_{i}\in S$ and prove
that there exists an optimal dividend payment strategy and has this form.
Later we extend the concept of strategies with this type of structure to the
case $S=[\underline{c},\overline{c}]$ by means of a curve in the state space
$[0,\infty)\times\lbrack\underline{c},\overline{c}]$ and look for the curve
which maximizes the expected discounted cumulative dividends.

\subsection{Optimal strategies for finite sets \label{Seccion optima Discreta}%
}

Take $S=\left\{  c_{1},c_{2},....,c_{n}\right\}  $ with $0\leq c_{1}%
<c_{2}<....<c_{n}=\overline{c}.$ Since for $i\leq n-1,$ the optimal value
function $V^{c_{i}}$ is a viscosity solution of (\ref{HJB equation discreta}),
there are values of $x$ where $V^{c_{i}}(x)=V^{c_{i+1}}(x)$ and values of $x$
where $\mathcal{L}^{c_{i}}(V^{c_{i}})(x)=0.$ We look for the simplest dividend
payment strategies, those whose value functions are solutions of
$\mathcal{L}^{c_{i}}=0$ for $x\in\lbrack0,z(c_{i}))$ and $V^{c_{i}}%
=V^{c_{i+1}}$ for $x\in\lbrack z(c_{i}),\infty)$ with some $z(c_{i})\geq0$. We
will show in this subsection that the optimal value function comes from such
types of strategies. More precisely, take $\widetilde{S}=S\setminus\{c_{n}\}$
and a function $z:\widetilde{S}\rightarrow\lbrack0,\infty)$; we define a
\textit{threshold strategy} by backward recursion, it is a stationary strategy
(which depends on both the current surplus $x$ and the implemented dividend
rate $c_{i}\in S$)%

\begin{equation}
\mathbf{\pi}^{z}=(C_{x,c_{i}})_{(x,c_{i})\in\lbrack0,\infty)\times S}\text{
where }C_{x,c_{i}}\in\Pi_{x,c_{i}}^{S} \label{Pizzeta}%
\end{equation}
as follows:

\begin{itemize}
\item If $i=n$, pay dividends with rate $c_{n}=$ $\overline{c}$ up to the time
of ruin, that is $(C_{x,c_{n}})_{t}=\overline{c}$.

\item If $1\leq i<n$ and $x\geq z(c_{i})$ follow $C_{x,c_{i+1}}\in
\Pi_{x,c_{i+1}}^{S}.$

\item If $1\leq i<n$ and $x<z(c_{i})$ pay dividends with rate $c_{i}$ as long
as the surplus is less than $z(c_{i})$ up to the ruin time; if the current
surplus reaches $z(c_{i})$ before the time of ruin, follow $C_{x,c_{i+1}}
\in\Pi_{x,c_{i+1}}^{S}$. More precisely
\[
(C_{x,c_{i}})_{t}=c_{i}I_{t\leq\tau\wedge\widehat{\tau}}+(C_{X_{\widehat{\tau
}},c_{i+1}})_{t~}I_{\widehat{\tau}\leq t<\tau},
\]
where $\widehat{\tau}$ is the first time at which the surplus reaches
$z(c_{i}).$
\end{itemize}

Let us call the value $z(c_{i})$ the \textit{threshold at dividend rate level}
$c_{i}$ and $z:\widetilde{S}\rightarrow\lbrack0,\infty)$ the \textit{threshold
function. }The value function of the stationary strategy $\mathbf{\pi}^{z}$ is
defined as
\begin{equation}
W^{z}(x,c_{i}):=J(x;C_{x,c_{i}}). \label{Wz discreta}%
\end{equation}
Note that $W^{z}(x,c_{i})$ only depends on $z(c_{k})$ for $i\leq k<n$,
$W^{z}(0,c_{i})=0$ and $W^{z}(x,c_{i})=V^{c_{n}}(x)$ for $x\geq\max
\{z(c_{k}):i<k<n\}.$

\begin{proposition}
\label{Formula recursiva de Wz discreta} We have the following recursive
formula for $W^{z}$:
\begin{align*}
W^{z}(x,c_{n})  &  =\frac{c_{n}}{q}\left(  1-e^{\theta_{2}(c_{n})x}\right)
,\\
W^{z}(x,c_{i})  &  =\left\{
\begin{array}
[c]{lll}%
W^{z}(x,c_{i+1}) & \text{if} & x\geq z(c_{i})\\
\frac{c_{i}}{q}\left(  1-e^{\theta_{2}(c_{i})x}\right)  +a^{z}(c_{i}
)(e^{\theta_{1}(c_{i})x}-e^{\theta_{2}(c_{i})x})\  & \text{if} & x<z(c_{i})
\end{array}
\right.
\end{align*}
for $i<n$, where
\[
a^{z}(c_{i}):=\frac{W^{z}(z(c_{i}),c_{i+1})-\frac{c_{i}}{q}\left(
1-e^{\theta_{2}(c_{i})z(c_{i})}\right)  }{e^{\theta_{1}(c_{i})z(c_{i}
)}-e^{\theta_{2}(c_{i})z(c_{i})}}.
\]

\end{proposition}

\textit{Proof.} We have that $\mathcal{L}^{c_{i}}(W^{z})(x,c_{i})=0$ for
$x\in(0,z(c_{i}))$ because the stationary strategy $\mathbf{\pi}^{z}$ pays
$c_{i}$ when the current surplus is in $(0,z(c_{i}))$. Also $W^{z}(0,c_{i})=0$
because ruin is immediate at $x=0$, and by definition $W^{z}(z(c_{i}%
),c_{i})=W^{z}(z(c_{i}),c_{i+1}).$ From
(\ref{Solucion de L=0 con condicion en 0}), we get the result. \hfill
$\blacksquare$\newline

Let us now look for the maximum of the value functions $W^{z}(x,c_{i})$ among
all the possible threshold functions\textit{\ }$z:\widetilde{S}\rightarrow
\lbrack0,\infty)$, and denote by $z^{\ast}$ the optimal threshold function.
From Proposition \ref{Estrategia trivial}, $z^{\ast}=0$ for $\overline
{c}=c_{n}\leq q\sigma^{2}/(2\mu)$, so that from now on we only consider the
case $c_{n}>q\sigma^{2}/(2\mu).$

Since the function $W^{z}(x,c_{n})$ is known, there are two ways to solve this
optimization problem (using a backward recursion). We will study the problem
using both of them.

\begin{enumerate}
\item The first approach consists of seeing the optimization problem as a
sequence of $n-1$ one-dimensional optimization problems, that is obtaining the
maximum $a^{z}(c_{i})$ for $i=n-1,\ldots,1$. If $W^{z^{\ast}}(x,c_{k})$ and
$z^{\ast}(c_{k})$ are known for $k=i+1,\ldots,n$, then from Proposition
\ref{Formula recursiva de Wz discreta} we can obtain
\begin{equation}
z^{\ast}(c_{i})=\min\left(  \arg\max_{y\in\lbrack0,\infty)}\frac{W^{z^{\ast}
}(y,c_{i+1})-\frac{c_{i}}{q}\left(  1-e^{\theta_{2}(c_{i})y}\right)
}{e^{\theta_{1}(c_{i})y}-e^{\theta_{2}(c_{i})y}}\right)  .
\label{First Approach}%
\end{equation}
Note that
\[
\lim_{y\rightarrow\infty}\frac{W^{z^{\ast}}(y,c_{i+1})-\frac{c_{i}}{q}\left(
1-e^{\theta_{2}(c_{i})y}\right)  }{e^{\theta_{1}(c_{i})y}-e^{\theta_{2}
(c_{i})y}}=0
\]
$\ $because $\lim_{y\rightarrow\infty}W^{z^{\ast}}(y,c_{i+1})=\frac{c_{n}}{q}
$, so $z^{\ast}(c_{i})$ exists.

\item As a second approach, one can view the optimization problem as a
backward recursion of obstacle problems (see Remark \ref{remark del obstaculo}
). If $W^{z^{\ast}}(x,c_{k})$ and $z^{\ast}(c_{k})$ are known for
$k=i+1,\ldots,n$, we look for the smallest solution $U^{\ast}$ of the equation
$\mathcal{L}^{c_{i}}(U)=0$ in $[0,\infty)$ with boundary condition $U(0)=0$
above $W^{z^{\ast}}(\cdot,c_{i+1}).$ Then
\begin{equation}
z^{\ast}(c_{i})=\min\{y>0:U^{\ast}(y)=W^{z^{\ast}}(y,c_{i+1} )\}.
\label{Second Approach}%
\end{equation}
By (\ref{Solucion de L=0 con condicion en 0}), the solutions $U$ of the
equation $\mathcal{L}^{c_{i}}(U)=0$ in $[0,\infty)$ with boundary condition
$U(0)=0$ are of the form
\[
U_{a}(x)=\frac{c_{i}}{q}\left(  1-e^{\theta_{2}(c_{i})x}\right)
+a(e^{\theta_{1}(c_{i})x}-e^{\theta_{2}(c_{i})x}).
\]
Hence $U_{a}(x)$ is increasing in $a$ and $\lim_{a\rightarrow\infty}
U_{a}(x)=\infty$ for $x>0$, and so there exists an $a_{i}^{\ast}>0$ such that
\[
U^{\ast}=U_{a_{i}^{\ast}}=\min\{U_{a}:U_{a}(x)\geq W^{z^{\ast}}(x,c_{i+1}
)\text{ for all }x\geq0\},
\]
because $\lim_{x\rightarrow\infty}U_{0}(x)=\frac{c_{i}}{q}<\frac{c_{n}}
{q}=\lim_{x\rightarrow\infty}W^{z^{\ast}}(x,c_{i+1})$.
\end{enumerate}

\begin{remark}
\normalfont
\label{Remark Obstaculo} In the second approach, we can see $z^{\ast}(c_{i})$
as the smallest $x>0$ such that $U_{a^{\ast}}$ and $W^{z^{\ast}}(\cdot
,c_{i+1})$ coincide; more precisely $U_{a^{\ast}}(z^{\ast}(c_{i}))=W^{z^{\ast
}}(z^{\ast}(c_{i}),c_{i+1})$, $U_{a^{\ast}}(x)\geq W^{z^{\ast}}(x,c_{i+1})$
for $x>0$ and $U_{a^{\ast}}(x)>W^{z^{\ast}}(x,c_{i+1})$ for $x\in(0,z^{\ast
}(c_{i}))$. Note that $U_{a^{\ast}}^{\prime}(z^{\ast}(c_{i}))=\partial
_{x}W^{z^{\ast}}(z^{\ast}(c_{i}),c_{i+1})$ and $U_{a^{\ast}}(\cdot
)-W^{z^{\ast}}(\cdot,c_{i+1})\ $is locally convex at $z^{\ast}(c_{i})$. By the
recursive construction, this implies that $W^{z^{\ast}}(x,c_{i})$ is
infinitely continuously differentiable at all $x\in\lbrack0,\infty
)\setminus\{z^{\ast}(c_{k}):k=i,...,n-1\}$ and continuously differentiable at
the points $z^{\ast}(c_{k})\ $for $k=i,...,n-1.$
\end{remark}

\begin{lemma}
\label{Remark Concavity} $U_{0}(x)$ is an increasing concave function. If
$a>0$, $U_{a}(x)$ is increasing, and is concave in $(-\infty,y_{0})$ and
convex in $(y_{0},\infty)$ with
\[
y_{0}:=\frac{\log\left(  \left(  \frac{c}{q}+a\right)  \theta_{2}(c)^{2}
\right)  -\log(a\theta_{1}(c)^{2})}{\theta_{1}(c)-\theta_{2}(c)}.
\]
In the case $c\leq\mu,$ we have that $y_{0}>0$; in the case $c>\mu,$ we have
that $y_{0}\leq0$ if and only if
\[
0<a<\frac{c}{q}\frac{\theta_{2}(c)^{2}}{\left(  \theta_{1}(c)^{2}-\theta
_{2}(c)^{2}\right)  }\text{.}
\]

\end{lemma}

\textit{Proof.} We have that%
\[
\partial_{x}U_{a}(x)=-\left(  \frac{c}{q}+a\right)  \theta_{2}(c)e^{\theta
_{2}(c)x}+a\theta_{1}(c)e^{\theta_{1}(c)x}>0,
\]
and%

\[
\partial_{xx}U_{a}(x)=-\left(  \frac{c}{q}+a\right)  \theta_{2}(c)^{2}%
e^{\theta_{2}(c)x}+a\theta_{1}(c)^{2}e^{\theta_{1}(c)x}\geq0
\]
if and only if $x\geq y_{0}$. The result follows from Definition
\ref{Definicion tita1 tita2}.\hfill$\blacksquare$\newline

In the next theorem, we show that there exists an optimal strategy and it is
of threshold type.

\begin{theorem}
\label{Optima Discreta Threshold} If $z^{\ast}\ $is the optimal threshold
function, then $W^{z^{\ast}}(x,c_{i})$ is the optimal function $V^{c_{i}}(x)$
defined in (\ref{Optimal Value Function}) for $i=1,...,n$.
\end{theorem}

\textit{Proof.} By definition $W^{z^{\ast}}(\cdot,c_{n})=V^{c_{n}}.$ Assuming
that $W^{z^{\ast}}(\cdot,c_{i+1})=V^{c_{i+1}}$ for $i=n-1,...,1$, by Theorem
\ref{Caracterizacion Discreta}, it is enough to prove that $W^{z^{\ast}}%
(\cdot,c_{i})$ is a viscosity solution of (\ref{HJB equation discreta}). Since
by construction $V^{c_{i+1}}-W^{z^{\ast}}(\cdot,c_{i})\leq0$, it remains to be
seen that $\mathcal{L}^{c_{i}}(W^{z^{\ast}})(x,c_{i})\leq0$ for $x\geq
z^{\ast}(c_{i})$. By Remark \ref{Remark Obstaculo}, $W^{z^{\ast}}(\cdot
,c_{i})$ is continuously differentiable and it is piecewise infinitely
differentiable in open intervals in which it solves $\mathcal{L}^{c_{j}%
}(W^{z^{\ast}})(x,c_{i})=0$ for some $j\geq i.$ By the definition of a
viscosity solution, it is enough to prove the result in these open intervals.
For $x$ in these open intervals,%

\[
\mathcal{L}^{c_{i}}(W^{z^{\ast}})(x,c_{i})=\mathcal{L}^{c_{j}}(W^{z^{\ast}%
})(x,c_{i})+(c_{i}-c_{j})(1-\partial_{x}W^{z^{\ast}}(x,c_{i}))\leq0
\]
if and only if $\partial_{x}W^{z^{\ast}}(x,c_{i})\leq1.$ There exists
$\delta>0$ and some $j>i$ such that $\mathcal{L}^{c_{j}}(W^{z^{\ast}}%
)(x,c_{i})=0$ in $(z^{\ast}(c_{i}),z^{\ast}(c_{i})+\delta)$ and then%

\[
\mathcal{L}^{c_{j}}(W^{z^{\ast}})(z^{\ast}(c_{i})^{+},c_{i})=0\text{,
}\mathcal{L}^{c_{i}}(W^{z^{\ast}})(z^{\ast}(c_{i})^{-},c_{i})=0,
\]
so%
\[%
\begin{array}
[c]{lll}%
0 & = & \mathcal{L}^{c_{i}}(W^{z^{\ast}})(z^{\ast}(c_{i})^{-},c_{i}
)-\mathcal{L}^{c_{j}}(W^{z^{\ast}})(z^{\ast}(c_{i})^{+},c_{i})\\
& = & \frac{\sigma^{2}}{2}(\partial_{xx}W^{z^{\ast}}(z^{\ast}(c_{i})^{-}
,c_{i})-\partial_{xx}W^{z^{\ast}}(z^{\ast}(c_{i})^{+-},c_{i}))\\
&  & +(c_{i}-c_{j})(1-\partial_{x}W^{z^{\ast}}(z^{\ast}(c_{i}),c_{i})).
\end{array}
\]
By Remark \ref{Remark Obstaculo}, $\partial_{xx}W^{z^{\ast}}(z^{\ast}%
(c_{i})^{-},c_{i})-\partial_{xx}W^{z^{\ast}}(z^{\ast}(c_{i})^{+-},c_{i})\geq0$
and $c_{i}-c_{j}<0,$ so we conclude that $\partial_{x}W^{z^{\ast}}(z^{\ast
}(c_{i}),c_{i}))\leq1$.

If $i=n-1$, $W^{z^{\ast}}(\cdot,c_{n-1})=W^{z^{\ast}}(\cdot,c_{n})$ for $x\geq
z^{\ast}(c_{n-1})$, by Remark \ref{Remark Concavity}, $W^{z^{\ast}}%
(\cdot,c_{n})$ is concave and so $W^{z^{\ast}}(x,c_{n})\leq W^{z^{\ast}%
}(z^{\ast}(c_{n-1}),c_{n})\leq1$ and we have the result.\newline

We need to prove that $\partial_{x}W^{z^{\ast}}(x,c_{i})=\partial
_{x}W^{z^{\ast}}(x,c_{i+1})\leq1$ for $x\geq z^{\ast}(c_{i}).$ By induction
hypothesis, we know that $\partial_{x}W^{z^{\ast}}(x,c_{i+1})=\partial
_{x}V^{c_{i+1}}\leq1$ for $x\geq z^{\ast}(c_{i+1})$. In the case that
$z^{\ast}(c_{i})\geq z^{\ast}(c_{i+1}),$ it is straightforward; in the case
that $z^{\ast}(c_{i})<z^{\ast}(c_{i+1})$, it is enough to prove it in the
interval $(z^{\ast}(c_{i}),z^{\ast}(c_{i+1}))$. But $\partial_{x}W^{z^{\ast}%
}(z^{\ast}(c_{i}),c_{i})\leq1$, $\partial_{x}W^{z^{\ast}}(z^{\ast}%
(c_{i+1}),c_{i})=\partial_{x}W^{z^{\ast}}(z^{\ast}(c_{i+1}),c_{i+1})\leq1,$
and by Lemma \ref{Remark Concavity} \ $\partial_{x}W^{z^{\ast}}(x,c_{i})$ is
either increasing, or decreasing, or decreasing and then increasing in the
interval $(z^{\ast}(c_{i}),z^{\ast}(c_{i+1}))$, so that we have the result.
\hfill$\blacksquare$\newline

Taking the derivative in (\ref{First Approach}) with respect to $y$, we get
implicit equations for the optimal threshold strategy.

\begin{proposition}
\label{Formula threshold discreta} $z^{\ast}(c_{i})\ $satisfies the implicit
equation
\[%
\begin{array}
[c]{c}%
0=\frac{c_{i}}{q}\theta_{2}(c_{i})e^{\theta_{2}(c_{i})y}(e^{\theta_{1}
(c_{i})y}-e^{\theta_{2}(c_{i})y})-\frac{c_{i}}{q}\left(  1-e^{\theta_{2}
(c_{i})y}\right)  \left(  \theta_{1}(c_{i})e^{\theta_{1}(c_{i})y}-\theta
_{2}(c_{i})e^{\theta_{2}(c_{i})y}\right) \\
+\partial_{x}W^{z^{\ast}}(y,c_{i+1})\left(  e^{\theta_{1}(c_{i})y}
-e^{\theta_{2}(c_{i})y}\right)  -W^{z^{\ast}}(y,c_{i+1})\left(  \theta
_{1}(c_{i})e^{\theta_{1}(c_{i})y}-\theta_{2}(c_{i})e^{\theta_{2}(c_{i}
)y}\right)
\end{array}
\]
for $i=n-1,\ldots,1.$
\end{proposition}

\begin{remark}
\normalfont
\label{Pizzeta extendido} Given $z:\widetilde{S}\rightarrow\lbrack0,\infty)$,
we have defined in (\ref{Pizzeta}) a threshold strategy $\mathbf{\pi}
^{z}=(C_{x,c_{i}})_{(x,c_{i})\in\lbrack0,\infty)\times S}$, where $C_{x,c_{i}
}\in\Pi_{x,c_{i}}^{S}$ for $i=1,\ldots,n$. We can extend this threshold
strategy to
\begin{equation}
\widetilde{\mathbf{\pi}}^{z}=(C_{x,c})_{(x,c)\in\lbrack0,\infty)\times\lbrack
c_{1},c_{n}]}~\text{where~}C_{x,c}\in\Pi_{x,c}^{S} \label{Pizzeta Rulo}%
\end{equation}
as follows:

\begin{itemize}
\item If $c\in(c_{i},c_{i+1})$ and $x<z(c_{i})$, pay dividends with rate $c$
while the current surplus is less than $z(c_{i})$ up to the time of ruin. If
the current surplus reaches $z(c_{i})$ before the time of ruin, follow
$C_{z(c_{i} ),c_{i+1}}\in\Pi_{x,c_{i+1}}^{S}$.

\item If $c\in(c_{i},c_{i+1})$ and $x\geq z(c_{i})$ for $1\leq i<n$, follow
$C_{x,c_{i+1}}\in\Pi_{x,c_{i+1}}^{S}.$
\end{itemize}

The value function of the stationary strategy $\widetilde{\mathbf{\pi}}^{z}$
is defined as
\begin{equation}
J^{\widetilde{\mathbf{\pi}}^{z}}(x,c):=J(x;C_{x,c}):[0,\infty)\times\lbrack
c_{1},c_{n}]\rightarrow{\mathbb{R}.} \label{W Pizzeta Rulo}%
\end{equation}

\end{remark}

\subsection{Curve strategies and the optimal curve strategy
\label{Curve strategies and optimal curve}}

As it is typical for these type of problems, the way in which the optimal
value function $V(x,c)$ solves the HJB equation (\ref{HJB equation}) suggests
that the state space $[0,\infty)\times\lbrack\underline{c},\overline{c}]$ is
partitioned into two regions: a \textit{non-change dividend region
}$\mathcal{NC}^{\ast}$ in which the dividends are paid at constant rate and a
\textit{change dividend region }$\mathcal{CH}^{\ast}$ in which the rate of
dividends increases. Roughly speaking, the region $\mathcal{NC}^{\ast}$
consists of the points in the state space where $\mathcal{L}^{c}(V)=0$ and
$\partial_{c}V<0$ and $\mathcal{CH}^{\ast}$ consists of the points where
$\partial_{c}V=0$. We introduce a family of stationary strategies (or limit of
stationary strategies) where the change and non-change dividend payment
regions are connected and split by a free boundary curve. This family of
strategies is the analogue to the threshold strategies for finite $S$
introduced in Section \ref{Seccion optima Discreta}.

Let us consider the set
\begin{equation}
\mathcal{B}=\{\zeta~s.t.\text{~}\zeta:[\underline{c},\overline{c}
)\rightarrow\lbrack0,\infty)\text{ is Riemann integrable and c\`{a}dl\`{a}g
with }\lim_{c\rightarrow\overline{c}^{-}}\zeta(c)<\infty\}. \label{Conjunto B}%
\end{equation}
In the first part of this subsection, we define the $\zeta$\textit{-value
function} $W^{\zeta}$ associated to a curve
\[
\mathcal{R}(\zeta)=\left\{  (\zeta(c),c):c\in\lbrack\underline{c},\overline
{c})\right\}  \subset\lbrack0,\infty)\times\lbrack\underline{c},\overline{c})
\]
for $\zeta\in$ $\mathcal{B}$, and we will see that, in some sense, $W^{\zeta
}(x,c)$ is a (limit) value function of the strategy which pays dividends at
constant rate in the case that $x<\zeta(c)$ and otherwise increases the rate
of dividends. So the curve $\mathcal{R}(\zeta)$ splits the state space
$[0,\infty)\times\lbrack\underline{c},\overline{c}]$ into two connected
regions: $\mathcal{NC(\zeta)=}\{(x,c)\in\lbrack0,\infty)\times\lbrack
\underline{c},\overline{c}]:x<\zeta(c)\}$ where dividends are paid with
constant rate, and $\mathcal{CH(\zeta)=}\{(x,c)\in\lbrack0,\infty
)\times\lbrack\underline{c},\overline{c}]:x\geq\zeta(c)\}$ where the dividend
rate increases. In the second part of the subsection we then will look for the
$\zeta_{0}\in$ $\mathcal{B}$ that maximizes the $\zeta$\textit{-value
function} $W^{\zeta}$, using calculus of variations.\newline

Let us consider the following auxiliary functions $b_{0}$, $b_{1}%
:(0,\infty)\times\lbrack\underline{c},\overline{c}]\rightarrow{\mathbb{R}} $
\begin{equation}%
\begin{array}
[c]{ccc}%
b_{0}(x,c) & := & \frac{-\frac{1}{q}\left(  1-e^{\theta_{2}(c)x}\right)
+\frac{c}{q}\theta_{2}^{\prime}(c)e^{\theta_{2}(c)x}x}{e^{\theta_{1}
(c)x}-e^{\theta_{2}(c)x}},\\
b_{1}(x,c) & := & \frac{\left(  -\theta_{1}^{\prime}(c)e^{\theta_{1}
(c)x}+\theta_{2}^{\prime}(c)e^{\theta_{2}(c)x}\right)  x}{e^{\theta_{1}
(c)x}-e^{\theta_{2}(c)x}}.
\end{array}
\label{b0 y b1}%
\end{equation}
Both $b_{0}(x,c)$ and $b_{1}(x,c)$ are not defined in $x=0,$ so we extend the
definition as
\[
b_{0}(0,c)=\lim_{x\rightarrow0^{+}}b_{0}(x,c):=\frac{c\theta_{2}^{\prime
}(c)-\theta_{2}(c)}{q\left(  \theta_{1}(c)-\theta_{2}(c)\right)  }%
\]
and%
\[
b_{1}(0,c)=\lim_{x\rightarrow0^{+}}b_{1}(x,c):=\frac{\theta_{2}^{\prime
}(c)-\theta_{1}^{\prime}(c)}{\theta_{1}(c)-\theta_{2}(c)}.
\]

In order to define the $\zeta$-value function in the non-change region
$\mathcal{NC(\zeta)}$, we will define and study in the next technical lemma
the functions $H^{\zeta}$ and $A^{\zeta}$ for any $\zeta\in\mathcal{B}$.

\begin{lemma}
\label{Funcion de valor de curva} Given $\zeta\in\mathcal{B}$, the unique
continuous function $H^{\zeta}:\{(x,c)\in\lbrack0,\infty)\times\lbrack
\underline{c},\overline{c}]:x\leq\zeta(c)\}\rightarrow\lbrack0,\infty)$ which
satisfies for any $c\in\lbrack\underline{c},\overline{c})$ that
\[
\mathcal{L}^{c}(H^{\zeta})(x,c)=0\text{ for }0\leq x\leq\zeta(c)~
\]
with boundary conditions $H^{\zeta}(0,c)=0$, $H^{\zeta}(x,\overline
{c})=V^{\left\{  \overline{c}\right\}  }(x,\overline{c})$ and $\partial
_{c}H^{\zeta}(\zeta(c),c)=0\ $ at the points of continuity of $\zeta$ is given
by
\begin{equation}
H^{\zeta}(x,c)=\frac{c}{q}\left(  1-e^{\theta_{2}(c)x}\right)  +A^{\zeta
}(c)(e^{\theta_{1}(c)x}-e^{\theta_{2}(c)x}), \label{Definicion de Hz}%
\end{equation}
where
\begin{equation}
A^{\zeta}(c)=-\int_{c}^{\overline{c}}e^{-\int_{c}^{t}b_{1}(\zeta(s),s)ds}
b_{0}(\zeta(t),t)dt . \label{Definicion A(z)}%
\end{equation}
Moreover, $A^{\zeta}$ satisfies $A^{\zeta}(\overline{c})=0$, is differentiable
and satisfies
\begin{equation}
\left(  A^{\zeta}\right)  ^{\prime}(c)=b_{0}(\zeta(c),c)+b_{1}(\zeta
(c),c)A^{\zeta}(c), \label{Ecuacion Diferencial de Az}%
\end{equation}
at the points where $\zeta$ is continuous.
\end{lemma}

\textit{Proof.} Since $\mathcal{L}^{c}(H^{\zeta}(x,c))=0$ and $H^{\zeta
}(0,c)=0$, we can write by (\ref{Solucion de L=0 con condicion en 0})
\[
H^{\zeta}(x,c)=\frac{c}{q}\left(  1-e^{\theta_{2}(c)x}\right)  +A^{\zeta
}(c)(e^{\theta_{1}(c)x}-e^{\theta_{2}(c)x}),
\]
where $A^{\zeta}(c)$ should be defined in such a way that $A^{\zeta}%
(\overline{c})=0$ (because $H^{\zeta}(x,\overline{c})=V^{\left\{  \overline
{c}\right\}  }(x,\overline{c}))$ and
\[%
\begin{array}
[c]{lll}%
0=\partial_{c}H^{\zeta}(\zeta(c),c) & = & \frac{1}{q}\left(  1-e^{\theta
_{2}(c)\zeta(c)}\right)  -\frac{c}{q}\theta_{2}^{\prime}(c)e^{\theta
_{2}(c)\zeta(c)}\zeta(c)+A^{\zeta}(c)^{\prime}(e^{\theta_{1}(c)\zeta
(c)}-e^{\theta_{2}(c)\zeta(c)})\\
&  & +A^{\zeta}(c)(\theta_{1}^{\prime}(c)e^{\theta_{1}(c)\zeta(c)}-\theta
_{2}^{\prime}(c)e^{\theta_{2}(c)\zeta(c)})\zeta(c)
\end{array}
\]
at the points of continuity of $\zeta$. Hence,%
\[%
\begin{array}
[c]{lll}%
\left(  A^{\zeta}\right)  ^{\prime}(c) & = & \frac{-\frac{1}{q}\left(
1-e^{\theta_{2} (c)x}\right)  +\frac{c}{q}\theta_{2}^{\prime}(c)e^{\theta
_{2}(c)x}x} {e^{\theta_{1}(c)x}-e^{\theta_{2}(c)x}}+A^{\zeta}(c)\frac{\left(
-\theta_{1}^{\prime}(c)e^{\theta_{1}(c)x}+\theta_{2}^{\prime}(c)e^{\theta_{2}
(c)x}\right)  x}{e^{\theta_{1}(c)x}-e^{\theta_{2}(c)x}}(\ast)\\
& = & b_{0}(\zeta(c),c)+b_{1}(\zeta(c),c)A^{\zeta}(c).
\end{array}
\]
Solving this ODE with boundary condition $A^{\zeta}(\overline{c})=0$, we get
the result. \hfill$\blacksquare$\newline

Given $\zeta\in\mathcal{B}$, we define the $\zeta$\textit{-value function}
\begin{equation}
W^{\zeta}(x,c):=\left\{
\begin{array}
[c]{lll}%
H^{\zeta}(x,c)\text{ } & \text{if} & (x,c)\in\mathcal{NC(\zeta)},\\
H^{\zeta}(x,C(x,c)) & \text{if} & (x,c)\in\mathcal{CH(\zeta)},
\end{array}
\right.  \label{Definicion de Wz}%
\end{equation}
where $H^{\zeta}$ is defined in Lemma \ref{Funcion de valor de curva} and%
\begin{equation}
C(x,c):=\max\{h\in\lbrack c,\overline{c}]:\zeta(d)\leq x{\large \ }
\text{{\large for}}{\large \ }d\in\lbrack c,h)\} \label{Definicion de c(x)}%
\end{equation}
in the case that $x\geq\zeta(c)$ and $c\in\lbrack\underline{c},\overline{c}%
)$.\newline

In the next propositions we will show that the $\zeta$-value function
$W^{\zeta}$ is the value function of an extended threshold strategy in the
case that $\zeta$ is a step function, and the limit of value functions of
extended threshold strategies in the case that $\zeta\in\mathcal{B}$.

\begin{proposition}
\label{Extended threshold strategy} Given $z:\widetilde{S}\rightarrow
\lbrack0,\infty)$ and the corresponding extended threshold strategy
$\widetilde{\mathbf{\pi}}^{z}$ defined in Remark \ref{Pizzeta extendido}, let
us consider the associated step function $\zeta\in\mathcal{B}$ defined as
\[
\zeta(c):=
{\displaystyle\sum\limits_{i=1}^{n-1}}
z(c_{i})I_{[c_{i},c_{i+1})}.
\]
Then the stationary value function of the extended threshold strategy
$\widetilde{\mathbf{\pi}}^{z}$ is given by
\[
J^{\widetilde{\mathbf{\pi}}^{z}}(x,c)=W^{\zeta}(x,c).
\]

\end{proposition}

\textit{Proof.} The stationary value function is continuous and satisfies
$\mathcal{L}^{c}(W^{\zeta})(x,c)=0$ for $0\leq x\leq\zeta(c)$, $W^{\zeta
}(0,c)=0$, $W^{\zeta}(x,\overline{c})=V^{\left\{  \overline{c}\right\}
}(x,\overline{c})$ and $\partial_{c}W^{\zeta}(\zeta(c),c)=0$ for
$c\notin\widetilde{S}$. Also the right-hand derivatives $\partial_{c}W^{\zeta
}(\zeta(c_{i}),c_{i}^{+})=0$ for $i=1,...,n-1.$ So, by Lemma
\ref{Funcion de valor de curva}, we obtain that $W^{\zeta}(x,c)=H^{\zeta
}(x,c)$ if $x<\zeta(c).$ If $x\geq\zeta(c),$ the result follows from the
definition of $\widetilde{\mathbf{\pi}}^{z}$. \hfill$\blacksquare$\newline

In the previous proposition we showed that in the case where $\zeta$ is the
associated step function of $z$, the stationary strategy $\widetilde
{\mathbf{\pi}}^{z}$ consists of increasing immediately the divided rate from
$c$ to $C(x,c)$ for $(x,c)\in\mathcal{CH(\zeta)}$, paying dividends at rate
$c$ until either reaching the curve $\mathcal{R(\zeta)}$ or ruin (whatever
comes first) for $(x,c)\in\mathcal{C(\zeta)}$, and paying dividends at rate
$\overline{c}$ until the time of ruin for $c=\overline{c}$.\newline

In the next proposition we show that for any $\zeta\in\mathcal{B}$, the
$\zeta$-value function $W^{\zeta}$ is the limit of value functions of extended
threshold strategies.

\begin{proposition}
Given $\zeta\in\mathcal{B}$, there exists a sequence of right-continuous step
functions $\zeta_{n}:[\underline{c},\overline{c})\rightarrow\lbrack0,\infty)$
such that $W^{\zeta_{n}}(x,c)$ converges uniformly to $W^{\zeta}(x,c)$.
\end{proposition}

\textit{Proof.} Since $\zeta$ is a Riemann integrable c\`{a}dl\`{a}g function,
we can approximate it uniformly by right-continuous step functions. Namely,
take a sequence of finite sets $\mathcal{S}^{k}=\{c_{1}^{k},c_{2}^{k}%
,\cdots,c_{n_{k}}^{k}\}$ with $\underline{c}=c_{1}^{k}<c_{2}^{k}%
<\cdots<c_{n_{k}}^{k}=\overline{c}$, and consider the right-continuous step
functions%
\[
\zeta_{k}(c)=%
{\displaystyle\sum\limits_{i=1}^{n_{k}-1}}
\zeta(c_{i}^{k})I_{[c_{i}^{k},c_{i+1}^{k})},
\]
such that $\delta(\mathcal{S}^{k})=\max_{i=1,\cdots,n_{k}-1}(c_{i+1}^{k}%
-c_{i}^{k})\rightarrow0$. We have that $\zeta_{k}\rightarrow\zeta$ uniformly,
and so both $A^{\zeta_{k}}(c)\rightarrow A^{\zeta}(c)$ and $W^{\zeta_{k}%
}(x,c)$ $\rightarrow$ $W^{\zeta}(x,c)$ uniformly.\hfill$\blacksquare$\newline

We now look for the maximum of $W^{\zeta}$ among $\zeta\in\mathcal{B}$. We
will show that if there exists a function $\zeta_{0}\in\mathcal{B}$ such that%
\begin{equation}
A^{\zeta_{0}}(\underline{c})=\max\{A^{\zeta}(\underline{c}):\zeta
\in\mathcal{B}\}, \label{Definicion z0}%
\end{equation}
then $W^{\zeta_{0}}(x,c)\geq W^{\zeta}(x,c)$ for all $(x,c)\in\lbrack
0,\infty)\times\lbrack c,\overline{c})\ $and $\zeta\in\mathcal{B}.$ This
follows from (\ref{Definicion de Hz}) and the next lemma, in which we prove
that the function $\zeta_{0}$ which maximizes (\ref{Definicion z0}) also
maximizes $A^{\zeta}(c)$ for any $c\in\lbrack\underline{c},\overline{c}).$

\begin{lemma}
For a given $c\in\lbrack\underline{c},\overline{c})$, define
\[
\mathcal{B}_{c}=\{~\zeta~st.\text{~}\zeta:[c,\overline{c})\rightarrow
\lbrack0,\infty)\text{ is Riemann integrable and c\`{a}dl\`{a}g with }
\lim_{d\rightarrow\overline{c}^{-}}\zeta(d)<\infty\}.
\]
If $\zeta_{0}\in\mathcal{B}$ satisfies (\ref{Definicion z0}), then for any
$c\in\lbrack\underline{c},\overline{c})$
\[
A^{\zeta_{0}}(c)=\max\{A^{\zeta}(c):\zeta\in\mathcal{B}_{c}\}.\
\]

\end{lemma}

\textit{Proof.} Given $\zeta\in\mathcal{B},$ we can write
\[
A^{\zeta}(\underline{c})=\left(  -\int_{\underline{c}}^{c}e^{-\int
_{\underline{c}}^{t}b_{1}(\zeta(s),s)ds}b_{0}(\zeta(t),t)dt\right)  +\left(
e^{-\int_{\underline{c}}^{c}b_{1}(\zeta(s),s)ds}\right)  A^{\zeta}(c).
\]

So
\[
A^{\zeta_{0}}(\underline{c})=\left(  -\int_{\underline{c}}^{c}e^{-\int
_{\underline{c}}^{t}b_{1}(\zeta_{0}(s),s)ds}b_{0}(\zeta_{0}(t),t)dt\right)
+\left(  e^{-\int_{\underline{c}}^{c}b_{1}(\zeta_{0}(s),s)ds}\right)
\max_{\zeta\in\mathcal{B}_{c}}A^{\zeta}(c)\text{. }
\]
\hfill$\blacksquare$\newline

Assuming that $\zeta_{0}$ exists, we will use calculus of variations to obtain
an implicit equation for $A^{\zeta_{0}}$. First we prove the following
technical lemma.

\begin{lemma}
\label{Condicion 1 para z derivable} For any~$c\in\lbrack\underline
{c},\overline{c}],$ we have
\[
\partial_{x}b_{1}(x,c)<0~\text{for }x>0~\text{and }\partial_{x}b_{1}
(0^{+},c)<0.
\]

\end{lemma}

\textit{Proof.}%
\[%
\begin{array}
[c]{lll}%
\partial_{x}b_{1}(x,c) & = & \dfrac{\theta_{1}^{\prime}(c)e^{2\theta_{1}
(c)x}(-1+e^{-(\theta_{1}(c)-\theta_{2}(c))x}(1+(\theta_{1}(c)-\theta
_{2}(c))x))}{(e^{\theta_{1}(c)x}-e^{\theta_{2}(c)x})^{2}}\\
&  & +\dfrac{\theta_{2}^{\prime}(c)e^{2\theta_{2}(c)x}(-1+e^{(\theta
_{1}(c)-\theta_{2}(c))x}(1-(\theta_{1}(c)-\theta_{2}(c))x))}{(e^{\theta
_{1}(c)x}-e^{\theta_{2}(c)x})^{2}}\\
& < & 0
\end{array}
\]

and, by Remark \ref{Propiedades de Titas},
\[
\lim_{x\rightarrow0}\partial_{x}b_{1}(x,c)=-\frac{\theta_{1}^{\prime
}(c)+\theta_{2}^{\prime}(c)}{2}<0\text{. }
\]
\hfill$\blacksquare$\newline

Let us now find the implicit equation for $A^{\zeta_{0}}.$

\begin{proposition}
\label{Condicion de Optimo} If the function $\zeta_{0}$ defined in
(\ref{Definicion z0}) exists, then $A^{\zeta_{0}}(c)$ satisfies
\[
A^{\zeta_{0}}(c)=-\frac{\partial_{x}b_{0}(\zeta_{0}(c),c)}{\partial_{x}
b_{1}(\zeta_{0}(c),c)}
\]
for all $c\in\lbrack\underline{c},\overline{c})$. Moreover, $A^{\zeta_{0}
}(\overline{c})=0$ and $A^{\zeta_{0}}(c)>0$ for $c\in\lbrack\underline
{c},\overline{c})$.
\end{proposition}

\textit{Proof.} Consider any function $\zeta_{1}\in\mathcal{B}$ with
$\zeta_{1}(\overline{c})=0$ then%

\[%
\begin{array}
[c]{ccc}%
A^{\zeta_{0}+\varepsilon\zeta_{1}}(\underline{c}) & = & -\int_{\underline{c}
}^{\overline{c}}e^{-\int_{\underline{c}}^{c}b_{1}(\zeta_{0}(s)+\varepsilon
\zeta_{1}(s),s)ds}b_{0}(\zeta_{0}(c)+\varepsilon\zeta_{1}(c),c)\,dc.
\end{array}
\]
Taking the derivative with respect to $\varepsilon$ and taking $\varepsilon
=0$, we get%

\[%
\begin{array}
[c]{lll}%
0=\left.  \partial_{\varepsilon}\left(  A^{\zeta_{0}+\varepsilon\zeta_{1}
}\right)  (\underline{c})\right\vert _{\varepsilon=0} & = & \int
_{\underline{c}}^{\overline{c}}\left(  (\int_{\underline{c}}^{c}\partial
_{x}b_{1}(\zeta_{0}(s),s)\zeta_{1}(s)ds)e^{-\int_{\underline{c}}^{c}
b_{1}(\zeta_{0}(s),s)ds}b_{0}(\zeta_{0}(c),c)\right)  dc\\
&  & -\int_{\underline{c}}^{\overline{c}}\left(  e^{-\int_{\underline{c}}
^{c}b_{1}(\zeta_{0}(s),s)ds}\partial_{x}b_{0}(\zeta_{0}(c),c)\zeta
_{1}(c)\right)  dc.\\
& = & \int_{\underline{c}}^{\overline{c}}\left(  \partial_{x}b_{1}(\zeta
_{0}(c),c)\zeta_{1}(c)\left(  \int_{c}^{\overline{c}}e^{-\int_{\underline{c}
}^{u}b_{1}(\zeta_{0}(s),s)ds}b_{0}(\zeta_{0}(u),u)du\right)  \right)  dc\\
&  & -\int_{\underline{c}}^{\overline{c}}\left(  e^{-\int_{\underline{c}}
^{c}b_{1}(\zeta_{0}(s),s)ds}\partial_{x}b_{0}(\zeta_{0}(c),c)\zeta
_{1}(c)\right)  dc.
\end{array}
\]
And so,%

\[
0=-\int_{\underline{c}}^{\overline{c}}\left(  e^{-\int_{\underline{c}}
^{c}b_{1}(\zeta_{0}(s),s)ds}\partial_{x}b_{0}(\zeta_{0}(c),c)-\partial
_{x}b_{1}(\zeta_{0}(c),c)(\int_{c}^{\overline{c}}e^{-\int_{\underline{c}}
^{u}b_{1}(\zeta_{0}(s),s)ds}b_{0}(\zeta_{0}(u),u)du)\right)  \zeta_{1}(c)dc.
\]
Since this holds for any $\zeta_{1}\in\mathcal{B}$ with $\zeta_{1}%
(\overline{c})=0$, we obtain that for any $c\in\lbrack c,\overline{c})$%

\[%
\begin{array}
[c]{lll}%
0 & = & \partial_{x}b_{0}(\zeta_{0}(c),c)-\partial_{x}b_{1}(\zeta
_{0}(c),c)\left(  \int_{c}^{\overline{c}}e^{-\int_{c}^{u}b_{1}(\zeta
_{0}(s),s)ds}b_{0}(\zeta_{0}(u),u)du\right) \\
& = & \partial_{x}b_{0}(\zeta_{0}(c),c)-\partial_{x}b_{1}(\zeta_{0}
(c),c)A^{\zeta_{0}}(c).
\end{array}
\]
Using Lemma \ref{Condicion 1 para z derivable}, we get the implicit equation
for $\zeta_{0}$. By definition $A^{\zeta_{0}}(\overline{c})=0$. Now take
$c\in\lbrack\underline{c},\overline{c})$, and the constant step function
$\zeta\in\mathcal{B}$ defined as $\zeta\equiv x_{0}$ where $x_{0}$ satisfies%

\[
\frac{c}{q}\left(  1-e^{\theta_{2}(c)x_{0}}\right)  <\frac{\overline{c}}%
{q}\left(  1-e^{\theta_{2}(\overline{c})x_{0}}\right)  .
\]
Then
\[
A^{\zeta_{0}}(c)\geq A^{\zeta}(c)=\dfrac{\frac{\overline{c}}{q}\left(
1-e^{\theta_{2}(\overline{c})x_{0}}\right)  -\frac{c}{q}\left(  1-e^{\theta
_{2}(c)x_{0}}\right)  }{e^{\theta_{1}(c)x_{0}}-e^{\theta_{2}(c)x}}>0\text{. }
\]
\hfill$\blacksquare$\newline

\noindent From now on, we extend the definition of $\zeta_{0}$ to
$[\underline{c},\overline{c}]$ as
\[
\zeta_{0}(\overline{c}):=\lim_{d\rightarrow\overline{c}^{-}}\zeta_{0}(d).
\]
Since $A^{\zeta_{0}}(\overline{c})=0$, we get from Proposition
\ref{Condicion de Optimo}%
\begin{equation}
\partial_{x}b_{0}(\zeta(\overline{c}),\overline{c} )=0,
\label{Condicion de Optimo en cbarra}%
\end{equation}
and since $A^{\zeta_{0}}(c)>0$ for $c\in\lbrack\underline{c},\overline{c})$,
we obtain that
\[
\partial_{x}b_{0}(\zeta_{0}(c),c)>0\text{.}%
\]
In the next proposition we show that, under some assumptions, the function
$\zeta_{0}:$ $[\underline{c},\overline{c}]\rightarrow\lbrack0,\infty)$ is the
unique solution of the first order differential equation
\begin{equation}
\zeta^{\prime}(c)=\left(  \dfrac{-b_{0}~(\partial_{x}b_{1})^{2}+b_{1}
~\partial_{x}b_{0}~\partial_{x}b_{1}-\partial_{xc}b_{0}~\partial_{x}
b_{1}+\partial_{xc}b_{1}~\partial_{x}b_{0}}{\partial_{xx}b_{0}~\partial
_{x}b_{1}-\partial_{xx}b_{1}~\partial_{x}b_{0}}\right)  (\zeta(c),c)
\label{Ecuacion diferencial de z0}%
\end{equation}
with boundary condition (\ref{Condicion de Optimo en cbarra}).

\begin{proposition}
\label{condiciones para zo derivable} If $\zeta_{0}(c)\ $defined in
(\ref{Definicion z0}) satisfies
\begin{equation}
\left(  \partial_{xx}b_{0}~\partial_{x}b_{1}-\partial_{xx}b_{1}~\partial
_{x}b_{0}\right)  (\zeta_{0}(c),c)\neq0\text{,}
\label{Condicion 2 para z derivable}%
\end{equation}
then $\zeta_{0}$ is infinitely differentiable and it is the unique solution of
(\ref{Ecuacion diferencial de z0}) with boundary condition
(\ref{Condicion de Optimo en cbarra}).
\end{proposition}

\textit{Proof.} From (\ref{Ecuacion Diferencial de Az}), we have%

\begin{equation}
\left(  -\frac{\partial_{x}b_{0}(\zeta_{0}(c),c)}{\partial_{x}b_{1}(\zeta
_{0}(c),c)}\right)  ^{\prime}=b_{0}(\zeta_{0}(c),c)+b_{1}(\zeta_{0}
(c),c)\left(  \frac{\partial_{x}b_{0}(\zeta_{0}(c),c)}{-\partial_{x}
b_{1}(\zeta_{0}(c),c)}\right)  . \label{Exp1}%
\end{equation}
By Assumption (\ref{Condicion 2 para z derivable}), the function%
\[
G(\zeta,c):=-\frac{\partial_{x}b_{0}(\zeta,c)}{\partial_{x}b_{1}(\zeta,c)}%
\]
satisfies
\[
\partial_{\zeta}G(\zeta_{0}(c),c)=-\left(  \frac{\partial_{xx}b_{0}(\zeta
_{0}(c),c)\partial_{x}b_{1}(\zeta_{0}(c),c)-\partial_{xx}b_{1}(\zeta
_{0}(c),c)\partial_{x}b_{0}(\zeta_{0}(c),c)}{\partial_{x}b_{1}(\zeta
_{0}(c),c)^{2}}\right)  \neq0.
\]
Hence, from (\ref{Exp1}), we have that $\zeta_{0}(c)$ is differentiable and we
get the differential equation (\ref{Ecuacion diferencial de z0}) for
$\zeta_{0}$. We obtain by a recursive argument that $\zeta_{0}(c)$ is
infinitely differentiable.\hfill$\blacksquare$\newline

In the next proposition, we state that the value function $W^{\zeta_{0}}$
satisfies a smooth-pasting property on the smooth free-boundary curve
\[
\mathcal{R(}\zeta_{0})=\{(\zeta_{0}(c),c)~\text{with }c\in\lbrack\underline
{c},\overline{c})\}.
\]
We also show that, under some conditions, $\zeta_{0}$ is the unique continuous
function $\zeta\in\mathcal{B}$ such that the associated $\zeta$-value function
$W^{\zeta}$ satisfies the smooth-pasting property at the curve $\mathcal{R(}%
\zeta)$.

\begin{proposition}
If $\zeta_{0}$ defined in (\ref{Definicion z0}) satisfies
(\ref{Condicion 2 para z derivable}), then $W^{\zeta_{0}}$ satisfies the
smooth-pasting property
\[
W_{cx}^{\zeta_{0}}(\zeta_{0}(c),c)=W_{cc}^{\zeta_{0}}(\zeta_{0}(c),c)=0\text{
for }c\in\lbrack\underline{c},\overline{c}]\text{.}
\]
Conversely, let $h:[0,\infty)\times\lbrack\underline{c},\overline
{c}]\rightarrow\lbrack0,\infty)\ $with $h(x,\overline{c})=V^{\{\overline{c}
\}}(x,\overline{c})$ and $h(0,c)=0$ for $c\in\lbrack\underline{c},\overline
{c})$. Assume that for $c\in\lbrack\underline{c},\overline{c}),$
\[
\zeta(c):=\sup\left\{  y:\mathcal{L}^{c}(h)(x,c)=0\text{ for }0\leq x\leq
y\right\}  \text{ }
\]
is a positive and continuous function in $\mathcal{B}$ satisfying
\[
\partial_{c}h(\zeta(c),c)=\partial_{cx}h(\zeta(c),c)=0
\]
and $\left(  \partial_{xx}b_{0}~\partial_{x}b_{1}-\partial_{xx}b_{1}
~\partial_{x}b_{0}\right)  (\zeta(c),c)\neq0$ for~$c\in\lbrack\underline
{c},\overline{c})$; then $\zeta$ coincides with $\zeta_{0}$ and $h(x,c)=$
$W^{\zeta_{0}}(x,c)$ for $0\leq x\leq\zeta(c)$ and $c\in\lbrack\underline
{c},\overline{c}]$.
\end{proposition}

\textit{Proof.} Let us define for $x\geq0$ and $c\in\lbrack\underline
{c},\overline{c}]$ the function%

\begin{equation}
H(x,c):=\frac{c}{q}\left(  1-e^{\theta_{2}(c)x}\right)  +a(c)(e^{\theta
_{1}(c)x}-e^{\theta_{2}(c)x})\text{,} \label{H en SP}%
\end{equation}
where $a:[\underline{c},\overline{c}]\rightarrow\lbrack0,\infty)$ is a
function with $a(\overline{c})=0$. Note that $H$ satisfies $\mathcal{L}%
^{c}(W(x,c))=0$ for all $x\geq0,$ $H(0,c)=0$ and $H(x,\overline{c}%
)=\frac{\overline{c}}{q}\left(  1-e^{\theta_{2}(c)x}\right)  $. We have,%

\[
\partial_{x}H(x,c)=-\frac{c}{q}\theta_{2}(c)e^{\theta_{2}(c)x}+a(c)(\theta
_{1}(c)e^{\theta_{1}(c)x}-\theta_{2}(c)e^{\theta_{2}(c)x}).
\]
If $a(c)$ is differentiable,%

\begin{equation}
\partial_{c}H(x,c)=(e^{\theta_{1}(c)x}-e^{\theta_{2}(c)x})\left(
-b_{0}(x,c)-b_{1}(x,c)a(c)+a^{\prime}(c)\right)  , \label{Derivada c de H}%
\end{equation}
and%

\[%
\begin{array}
[c]{lll}%
\partial_{cx}H(x,c)=\partial_{xc}H(x,c) & = & (\theta_{1}(c)e^{\theta_{1}
(c)x}-\theta_{2}(c)e^{\theta_{2}(c)x})\left(  -b_{0}(x,c)-b_{1}
(x,c)a(c)+a^{\prime}(c)\right) \\
&  & +(e^{\theta_{1}(c)x}-e^{\theta_{2}(c)x})\left(  -\partial_{x}
b_{0}(x,c)-\partial_{x}b_{1}(x,c)a(c)\right)  .
\end{array}
\]
In the case that $a(c)=A^{\zeta_{0}}(c)$, take $H=H^{\zeta_{0}}$ as defined in
(\ref{Definicion de Hz}), by (\ref{Ecuacion Diferencial de Az}) and
Proposition \ref{Condicion de Optimo}, we obtain $\partial_{cx}H^{\zeta_{0}%
}(\zeta_{0}(c),c)=0$. Since $W^{\zeta_{0}}(x,c)=H^{\zeta_{0}}(x,c)$ for
$x<\zeta_{0}(c)$ and $W^{\zeta_{0}}(x,c)=H^{\zeta_{0}}(x,C(x,c))$ for
$x\geq\zeta_{0}(c)$, we get $\partial_{c}W^{\zeta_{0}}(x,c)=0$ for $x\geq
\zeta_{0}(c)$ and so $\partial_{cx}W^{\zeta_{0}}(\zeta(c),c)=0$. From
(\ref{Derivada c de H}), we get%

\[
\partial_{c}H^{\zeta_{0}}(x,c)=(e^{\theta_{1}(c)x}-e^{\theta_{2}(c)x})\left(
b_{0}(\zeta_{0}(c),c)-b_{0}(x,c)+\left(  b_{1}(\zeta_{0}(c),c)-b_{1}%
(x,c)\right)  A^{\zeta_{0}}(c)\right)  .
\]
Hence, from (\ref{Condicion 2 para z derivable}), we have that $\partial
_{cc}H^{\zeta_{0}}$ exists. Since $\partial_{c}H^{\zeta_{0}}(\zeta
_{0}(c),c)=0$ for $c\in\lbrack\underline{c},\overline{c}],$%
\[%
\begin{array}
[c]{lll}%
0 & = & \frac{d}{dc}(\partial_{c}H^{\zeta_{0}}(\zeta_{0}(c),c))\text{ }\\
& = & \partial_{cc}H^{\zeta_{0}}(\zeta_{0}(c),c)+\partial_{cx}H^{\zeta_{0}
}(\zeta_{0}(c),c)\zeta_{0}^{\prime}(c)\\
& = & \partial_{cc}H^{\zeta_{0}}(\zeta_{0}(c),c).
\end{array}
\]
Finally, since $W^{\zeta_{0}}(x,c)=H^{\zeta_{0}}(x,C(x,c))$ if $x\geq\zeta
_{0}(c)$, we get $\partial_{cc}W^{\zeta_{0}}(x,c)=0$ if $x\geq\zeta_{0}(c)$
and so $\partial_{cc}W^{\zeta_{0}}(\zeta_{0}(c),c)=0$.

Conversely, note that there exists $a(c)$ such that $h(x,c)=H(x,c)$ defined in
(\ref{H en SP}) for $x<\zeta(c)$; the existence of $\partial_{c}h$ implies
that $a(c)$ is differentiable. Hence,%

\[
0=\partial_{c}h(\zeta(c),c)=(e^{\theta_{1}(c)\zeta(c)}-e^{\theta_{2}
(c)\zeta(c)})\left(  -b_{0}(\zeta(c),c)-b_{1}(\zeta(c),c)a(c)+a^{\prime
}(c)\right)
\]
which implies%

\[
a^{\prime}(c)=b_{0}(\zeta(c),c)+b_{1}(\zeta(c),c)a(c).
\]
Also,%

\[%
\begin{array}
[c]{lll}%
0=\partial_{cx}h(\zeta(c),c) & = & (e^{\theta_{1}(c)\zeta(c)}-e^{\theta
_{2}(c)\zeta(c)})\left(  -\partial_{x}b_{0}(\zeta(c),c)-\partial_{x}
b_{1}(\zeta(c),c)a(c)\right) \\
&  & +(\theta_{1}(c)e^{\theta_{1}(c)\zeta(c)}-\theta_{2}(c)e^{\theta
_{2}(c)\zeta(c)})\left(  -b_{0}(\zeta(c),c)-b_{1}(\zeta(c),c)a(c)+a^{\prime
}(c)\right)
\end{array}
\]
implies%

\[
\partial_{x}b_{0}(\zeta(c),c)=\partial_{x}b_{1}(\zeta(c),c)a(c).
\]
Since $\left(  \partial_{xx}b_{0}~\partial_{x}b_{1}-\partial_{xx}%
b_{1}~\partial_{x}b_{0}\right)  (\zeta(c),c)\neq0$, both $\zeta$ and
$\zeta_{0}$ satisfy the same equation and so they coincide.\hfill
$\blacksquare$\newline

In the next proposition, we show more regularity for $W^{\zeta_{0}}$ in the
case that $\zeta_{0}$ is increasing.

\begin{proposition}
\label{z creciente} If $\zeta_{0}\ $defined in (\ref{Definicion z0}), is
increasing and satisfies (\ref{Condicion 2 para z derivable}), then
$W^{\zeta_{0}}$ is (2,1)-differentiable. Also, since the inverse $\zeta
_{0}^{-1}$ exists, $C(x,c)$ can be written in a simpler way:
\[
C(x,c)=\left\{
\begin{array}
[c]{lll}%
\overline{c} & \text{if} & \zeta_{0}(\overline{c})\leq x,\\
\zeta_{0}^{-1}(x) & \text{if} & \zeta_{0}(c)\leq x<\zeta_{0}(\overline{c}).
\end{array}
\right.
\]

\end{proposition}

\textit{Proof.} It is enough to prove that $\partial_{xx}W^{\zeta_{0}}%
(x^{+},c)=$ $\partial_{xx}W^{\zeta_{0}}(x^{-},c)$ for $\zeta_{0}(c)\leq
x<\zeta_{0}(\overline{c})$. We have%
\[%
\begin{array}
[c]{lll}%
\partial_{x}W^{\zeta_{0}}(x^{+},c) & = & \partial_{x}H^{\zeta_{0}}(x,\zeta
_{0}^{-1}(x))+\partial_{c}H^{\zeta_{0}}(x,\zeta_{0}^{-1}(x))\left(  \zeta
_{0}^{-1}\right)  ^{\prime}(x)\\
& = & \partial_{x}H^{\zeta_{0}}(x,\zeta_{0}^{-1}(x))\\
& = & \partial_{x}W^{\zeta_{0}}(x^{-},c).
\end{array}
\]
And so,
\[%
\begin{array}
[c]{lll}%
\partial_{xx}W^{\zeta_{0}}(x^{+},c) & = & \partial_{xx}H^{\zeta_{0}}
(x,\zeta_{0}^{-1}(x))+\partial_{cx}H^{\zeta_{0}}(x,\zeta_{0}^{-1}(x))\left(
\zeta_{0}^{-1}\right)  ^{\prime}(x)\\
& = & \partial_{xx}H^{\zeta_{0}}(x,\zeta_{0}^{-1}(x))\\
& = & \partial_{xx}W^{\zeta_{0}}(x^{-},c)\text{. }%
\end{array}
\]
\hfill$\blacksquare$

\subsection{Optimal strategies for the closed interval $S=[\underline
{c},\overline{c}]$ \label{Optimal strategies for the closed interval}}

First in this section, we give a verification result in order to check if a
$\zeta$-value function $W^{\zeta}$ is the optimal value function $V$. Our
conjecture is that the solution $\overline{\zeta}$ of
(\ref{Ecuacion diferencial de z0}) with boundary condition
(\ref{Condicion de Optimo en cbarra}) exists and is non-decreasing in
$[\underline{c},\overline{c}]$, and that $W^{\overline{\zeta}}$ coincides with
$V$, so that there exists an optimal curve strategy.

Using Proposition \ref{Estrategia trivial}, we know that the conjecture holds
for $\overline{c}\leq q\sigma^{2}/(2\mu)$ taking $\zeta_{0}\equiv0.$ In the
case that $\overline{c}>q\sigma^{2}/(2\mu)$, we will show that $\overline
{\zeta}$ exists and is increasing and $W^{\overline{\zeta}}=V$ for
$[\overline{c}-\varepsilon,\overline{c}]$ for $\varepsilon>0\ $small enough.
We were not able to prove the conjecture in the general case, although it
holds in our numerical explorations (see Section \ref{Numerical examples}).
However, we will prove that the $\zeta$-value functions $W^{\zeta}$ are
$\varepsilon$-optimal in the following sense: There exists a sequence
$\zeta_{n}\in\mathcal{B}$ such that $W^{\zeta_{n}}$ converges uniformly to the
optimal value function $V$.\newline

We state now a verification result for checking whether the $\zeta$-value
function $W^{\zeta}$ with $\zeta$ continuous is the optimal value function $V
$. In this verification result it is not necessary to use viscosity solutions
because the proposed value function solves the HJB equation in a classical
way. We will check these verification conditions for the limit value function
associated to the unique solution of (\ref{Ecuacion diferencial de z0}) with
boundary condition (\ref{Condicion de Optimo en cbarra}) (if it exists).

\begin{proposition}
\label{Verificacion Curva} If there exists a smooth function $\overline{\zeta
}$ such that the $\overline{\zeta}$-value function $W^{\overline{\zeta}}$ is
(2,1)-differentiable and satisfies
\[
\partial_{x}W^{\overline{\zeta}}(\overline{\zeta}(c),c)\leq1\text{ for }
c\in\lbrack\underline{c},\overline{c}]\quad\text{ and }\quad\partial
_{c}W^{\overline{\zeta}}(x,c)\leq0\text{ for }x\in\lbrack0,\overline{\zeta
}(c))\text{ and }c\in\lbrack\underline{c},\overline{c}),
\]
then $W^{\overline{\zeta}}=V$.
\end{proposition}

\textit{Proof.} We have that $\partial_{x}W^{\overline{\zeta}}(x,\overline
{c})\leq1$ for $x\geq\overline{\zeta}(\overline{c})$ because $W^{\overline
{\zeta}}(\cdot,\overline{c})$ is concave and $\partial_{x}W^{\overline{\zeta}%
}(\overline{\zeta}(\overline{c}),\overline{c})\leq1.$ Since $\mathcal{L}%
^{c}W^{\overline{\zeta}}(x,c)=0$ for $x\in\lbrack0,\overline{\zeta}(c)),$
$c\in\lbrack\underline{c},\overline{c})$; $\partial_{c}W^{\overline{\zeta}%
}(x,c)=0$ for $x\geq\overline{\zeta}(c)),$ $c\in\lbrack\underline{c}%
,\overline{c})$ and $W^{\overline{\zeta}}(\cdot,\overline{c})=V(\cdot
,\overline{c})$; by Theorem \ref{Caracterizacion Continua} it is sufficient to
prove that $\mathcal{L}^{c}W^{\overline{\zeta}}(x,c)\leq0$ for $x\geq
\overline{\zeta}(c)),$ $c\in\lbrack\underline{c},\overline{c})$. In this case,
we have that
\[
C(x,c)=\max\{h\in\lbrack c,\overline{c}]:\overline{\zeta}(d)\leq
x{\large \ }\text{{\large for}}{\large \ }d\in\lbrack c,h)\}
\]
satisfies $C(x,c)\geq c$, and also either $C(x,c)=\overline{c}$ or
$\overline{\zeta}(C(x,c))=x.$ So, we obtain $\mathcal{L}^{C(x,c)}V(x,C(x,c))=0
$ and then%

\[%
\begin{array}
[c]{lll}%
\mathcal{L}^{c}V(x,c) & = & \mathcal{L}^{C(x,c)}V(x,c)+(C(x,c)-c)(\partial
_{x}V(x,C(x,c))-1)\\
& = & (C(x,c)-c)(V_{x}(x,C(x,c))-1)\leq0\text{. }%
\end{array}
\]
\hfill$\blacksquare$\newline

Next we see that there exists a unique solution $\overline{\zeta}$ of
(\ref{Ecuacion diferencial de z0}) with boundary condition
(\ref{Condicion de Optimo en cbarra}) at least in $[\overline{c}%
-\varepsilon,\overline{c}]$ for some $\varepsilon>0$. First, let us study the
boundary condition (\ref{Condicion de Optimo en cbarra}) in the case
$\overline{c}$ $>q\sigma^{2}/(2\mu)$.

\begin{lemma}
\label{Punto inicial unico} If $\overline{c}$ $>q\sigma^{2}/(2\mu)$, there
exists a unique $\overline{z}>0$ such that $\partial_{x}b_{0} (\overline
{z},\overline{c})=0$; moreover $\partial_{x}b_{0}(x,\overline{c})<0$ for
$x\in\lbrack0,\overline{z})$ and $\partial_{x}b_{0}(x,\overline{c})>0$ for
$x\in(\overline{z},\infty)$. Also $\partial_{xx}b_{0}(\overline{z}
,\overline{c})>0.$
\end{lemma}

\textit{Proof.} For $x>0,$
\begin{subequations}
\begin{equation}
\partial_{x}b_{0}(x,c)=\tfrac{e^{\theta_{1}(c)x}\theta_{1}(c)-e^{\theta
_{2}(c)x}\theta_{2}(c)-ce^{2\theta_{2}(c)x}\theta_{2}^{\prime}(c)+e^{(\theta
_{1}(c)+\theta_{2}(c))x}(-\theta_{1}(c)+\theta_{2}(c)+c\theta_{2}^{\prime
}(c)(1-x\theta_{1}(c)+x\theta_{2}(c))}{q\left(  e^{\theta_{1}(c)x}
-e^{\theta_{2}(c)x}\right)  ^{2}}, \label{b0x}%
\end{equation}
and so
\end{subequations}
\[
\lim_{x\rightarrow0^{+}}\partial_{x}b_{0}(x,\overline{c})=-\frac{1}{2q}%
(\tfrac{\theta_{1}(\overline{c})\theta_{2}(\overline{c})}{\theta_{1}
(\overline{c})-\theta_{2}(\overline{c})}+\overline{c}\theta_{2}^{\prime
}(\overline{c}))=\tfrac{\overline{c}^{2}+q\sigma^{2}-\overline{c}(\mu
+\sqrt{(\overline{c}-\mu)^{2}+2q\sigma^{2}})}{2q\sigma^{2}\sqrt{(\overline
{c}-\mu)^{2}+2q\sigma^{2}}}.
\]
Hence, $\lim_{x\rightarrow0^{+}}\partial_{x}b_{0}(x,\overline{c})\geq0$ for
$\overline{c}\leq q\sigma^{2}/(2\mu)$ and $\lim_{x\rightarrow0^{+}}%
\partial_{x}b_{0}(x,\overline{c})<0$ for $\overline{c}>q\sigma^{2}/(2\mu)$.
Also
\[
\lim_{x\rightarrow\infty}\partial_{x}b_{0}(x,\overline{c})e^{\theta
_{1}(\overline{c})x}=\frac{\theta_{1}(\overline{c})}{q}>0.
\]
So, for $\overline{c}$ $>q\sigma^{2}/(2\mu)$ there exists (at least one)
$\overline{z}>0$ such that $\partial_{x}b_{0}(\overline{z},\overline{c})=0.$

We are showing next that $\partial_{x}b_{0}(x,\overline{c})=0$ for $x>0$
implies that $\partial_{xx}b_{0}(x,\overline{c})>0$. Consequently the result follows.

From (\ref{b0x}), we can write%
\[
\partial_{x}b_{0}(x,\overline{c})q\left(  e^{\theta_{1}(\overline{c}%
)x}-e^{\theta_{2}(\overline{c})x}\right)  ^{2}=g_{11}(x,\overline{c}%
)\theta_{2}^{\prime}(\overline{c})+g_{10}(x,\overline{c})
\]
and
\[
\partial_{xx}b_{0}(x,\overline{c})q\left(  e^{\theta_{1}(\overline{c}%
)x}-e^{\theta_{2}(\overline{c})x}\right)  ^{3}=g_{21}(x,\overline{c}%
)\theta_{2}^{\prime}(\overline{c})+g_{20}(x,\overline{c}),
\]
where%
\[
g_{11}(x,\overline{c})=-\overline{c}e^{2\theta_{2}(\overline{c})x}%
(1-e^{(\theta_{1}(\overline{c})-\theta_{2}(\overline{c}))x}(1-x(\theta
_{1}(\overline{c})-\theta_{2}(\overline{c}))),
\]%
\[
g_{10}(x,\overline{c})=-\theta_{1}(\overline{c})e^{\theta_{1}(\overline{c}%
)x}\left(  e^{\theta_{2}(\overline{c})x}-1\right)  +\theta_{2}(\overline
{c})e^{\theta_{2}(\overline{c})x}\left(  e^{\theta_{1}(\overline{c}%
)x}-1\right)  ,
\]%
\[
g_{21}(x,\overline{c})=\overline{c}e^{(\theta_{1}(\overline{c})+\theta
_{2}(\overline{c}))x}(\theta_{1}(\overline{c})-\theta_{2}(\overline
{c}))\left(  e^{\theta_{1}(\overline{c})x}(-2+x\theta_{1}(\overline
{c})-x\theta_{2}(\overline{c}))+e^{\theta_{2}(\overline{c})x}(2+x\theta
_{1}(\overline{c})-x\theta_{2}(\overline{c}))\right)  ,
\]%
\[%
\begin{array}
[c]{lll}%
g_{20}(x,\overline{c}) & = & e^{\theta_{1}(\overline{c})x}\theta_{1}%
^{2}(\overline{c})(e^{\theta_{2}(\overline{c})x}-1)(e^{\theta_{1}(\overline
{c})x}+e^{\theta_{2}(\overline{c})x})-2e^{(\theta_{1}(\overline{c})+\theta
_{2}(\overline{c}))x}\theta_{1}(\overline{c})\theta_{2}(\overline
{c})(-2+e^{\theta_{1}(\overline{c})x}+e^{\theta_{2}(\overline{c})x})\\
&  & +e^{\theta_{2}(\overline{c})x}\theta_{2}^{2}(\overline{c})(e^{\theta
_{1}(\overline{c})x}-1)(e^{\theta_{1}(\overline{c})x}+e^{\theta_{2}%
(\overline{c})x}).
\end{array}
\]
If $x>0$, take $u=x(\theta_{1}(\overline{c})-\theta_{2}(\overline{c}))>0$, we
can write
\[
-\frac{g_{11}(x,\overline{c})}{\overline{c}e^{2\theta_{2}(\overline{c})x}%
}=1-e^{u}(1-u)>0
\]
which implies that $g_{11}(x,\overline{c})<0.$

Consider now%
\[%
\begin{array}
[c]{lll}%
g(x,\overline{c}) & := & \partial_{xx}b_{0}(x,\overline{c})q\left(
e^{\theta_{1}(\overline{c})x}-e^{\theta_{2}(\overline{c})x}\right)  ^{3}
g_{11}(x,\overline{c})-\partial_{x}b_{0}(x,\overline{c})q\left(  e^{\theta
_{1}(\overline{c})x}-e^{\theta_{2}(\overline{c})x}\right)  ^{2}g_{21}
(x,\overline{c})\\
& = & g_{20}(x,\overline{c})g_{11}(x,\overline{c})-g_{10}(x,\overline
{c})g_{21}(x,\overline{c}).
\end{array}
\]
We are going to prove that $g(x,\overline{c})<0$ for $x>0$. For that purpose,
take%
\[
g_{0}(x,\overline{c}):=-\frac{g(x,\overline{c})}{\overline{c}\left(
e^{\theta_{1}(\overline{c})x}-e^{\theta_{2}(\overline{c})x}\right)
^{2}e^{2\theta_{2}(\overline{c})x}\theta_{1}^{2}(\overline{c})}\text{.}%
\]
Calling $t=-\frac{\theta_{2}(\overline{c})}{\theta_{1}(\overline{c})}>0$ and
$s=\theta_{1}(\overline{c})x>0$, and we can write%
\[
g_{0}(x,\overline{c})=-t^{2}+e^{s}(1+t)^{2}+e^{s+st}(-1-2t+st+st^{2}).
\]
Then $g_{0}(x,\overline{c})>0$ for $x>0$, because $g_{0}\left(  s,0\right)
=\partial_{t}g_{0}\left(  s,0\right)  =0$, $\partial_{t}^{2}g_{0}\left(
s,0\right)  =-2+(2-2s+s^{2})e^{s}>0$ and
\[
\partial_{t}^{3}g_{0}\left(  s,t\right)  =s^{3}e^{s+st}(2+4t+st+st^{2})>0
\]
for $s,t>0$. Finally, if $\partial_{x}b_{0}(\overline{z},\overline{c})=0$ for
$\overline{z}>0$, since
\[%
\begin{array}
[c]{lll}%
g(\overline{z},\overline{c}) & = & \partial_{xx}b_{0}(\overline{z}
,\overline{c})q\left(  e^{\theta_{1}(\overline{c})\overline{z}}-e^{\theta
_{2}(\overline{c})\overline{z}}\right)  ^{3}g_{11}(\overline{z},\overline
{c})-\partial_{x}b_{0}(\overline{z},\overline{c})q\left(  e^{\theta
_{1}(\overline{c})\overline{z}}-e^{\theta_{2}(\overline{c})\overline{z}
}\right)  ^{2}g_{21}(\overline{z},\overline{c})\\
& = & \partial_{xx}b_{0}(\overline{z},\overline{c})q\left(  e^{\theta
_{1}(\overline{c})\overline{z}}-e^{\theta_{2}(\overline{c})\overline{z}
}\right)  ^{3}g_{11}(\overline{z},\overline{c}),
\end{array}
\]
$g(\overline{z},\overline{c})<0$ and $g_{11}(\overline{z},\overline{c})$ $<0$,
we get that $\partial_{xx}b_{0}(\overline{z},\overline{c})>0$.\hfill
$\blacksquare$

\begin{proposition}
\label{Existencia cerca de Cbarra} In the case $\overline{c}$ $>q\sigma
^{2}/(2\mu)$ there exists a unique increasing solution $\overline{\zeta}$ of
(\ref{Ecuacion diferencial de z0}) with boundary condition
(\ref{Condicion de Optimo en cbarra}) in $[\overline{c}-\varepsilon
,\overline{c}] $ for some $\varepsilon>0$.
\end{proposition}

\textit{Proof.} From Lemma \ref{Punto inicial unico}, $\overline{\zeta
}(\overline{c})=\overline{z}$. By (\ref{Ecuacion diferencial de z0}) and since
the functions $b_{0}$ and $b_{1}$ are infinitely differentiable, it suffices
to prove that%

\[
\left(  \partial_{xx}b_{0}~\partial_{x}b_{1}-\partial_{xx}b_{1}~\partial
_{x}b_{0}\right)  (x,c)\neq0
\]
in a neighborhood of $\left(  \overline{z},\overline{c}\right)  .$ From Lemmas
\ref{Condicion 1 para z derivable} and \ref{Punto inicial unico},
\[
\left(  \partial_{xx}b_{0}~\partial_{x}b_{1}-\partial_{xx}b_{1}~\partial
_{x}b_{0}\right)  \left(  \overline{z},\overline{c}\right)  =(\partial
_{xx}b_{0}~\partial_{x}b_{1})\left(  \overline{z},\overline{c}\right)  <0.
\]
The existence of $\overline{\zeta}$ follows by continuity.

In order to show that $\overline{\zeta}$ is increasing in $[\overline
{c}-\varepsilon,\overline{c}]$ for some $\varepsilon>0$, it is sufficient to
prove
\[
\left(  -b_{0}~(\partial_{x}b_{1})^{2}+b_{1}~\partial_{x}b_{0}~\partial
_{x}b_{1}-\partial_{xc}b_{0}~\partial_{x}b_{1}+\partial_{xc}b_{1}~\partial
_{x}b_{0}\right)  (x,c)<0
\]
in a neighborhood of $\left(  \overline{z},\overline{c}\right)  .$ Since
$\partial_{x}b_{0}(\overline{z},\overline{c})=0$, we get%

\[%
\begin{array}
[c]{l}%
\left(  -b_{0}~(\partial_{x}b_{1})^{2}+b_{1}~\partial_{x}b_{0}~\partial
_{x}b_{1}-\partial_{xc}b_{0}~\partial_{x}b_{1}+\partial_{xc}b_{1}~\partial
_{x}b_{0}\right)  (\overline{z},\overline{c})\\
=-\partial_{x}b_{1}(\overline{z},\overline{c})(b_{0}(\overline{z},\overline
{c})\partial_{x}b_{1}(\overline{z},\overline{c})+\partial_{xc}b_{0}
(\overline{z},\overline{c})~);
\end{array}
\]
and since $\partial_{x}b_{1}(\overline{z},\overline{c})<0$, it is enough to
show that%
\[
b_{0}(\overline{z},\overline{c})\partial_{x}b_{1}(\overline{z},\overline
{c})+\partial_{xc}b_{0}(\overline{z},\overline{c})~<0.
\]
Taking $t=-\frac{\theta_{2}(\overline{c})}{\theta_{1}(\overline{c})}>0$ and
$u=\dfrac{-q~}{\theta_{2}(\overline{c})\sigma^{2}}\overline{z}>0$, we can
write%
\[
\frac{b_{0}(\overline{z},\overline{c})\partial_{x}b_{1}(\overline{z}
,\overline{c})+\partial_{xc}b_{0}(\overline{z},\overline{c})}{g_{0}
(\overline{z},\overline{c})}=g_{1}(u,t)+\frac{\overline{c}\overline{z}}%
{\sigma^{2}}g_{2}(u,t)
\]
and
\[
\frac{\partial_{x}b_{0}(\overline{z},\overline{c})}{f_{0}(\overline
{z},\overline{c})}=f_{1}(u,t)+\frac{\overline{c}\overline{z}}{\sigma^{2}}%
f_{2}(u,t)=0;
\]
where%
\[
g_{0}(\overline{z},\overline{c})=\frac{\overline{z}\left(  e^{\theta
_{1}(\overline{c})\overline{z}}-e^{\theta_{2}(\overline{c})\overline{z}
}\right)  ^{3}q^{2}(\theta_{1}(\overline{c})-\theta_{2}(\overline{c}))^{3}%
}{2\theta_{1}^{3}(\overline{c})(-\theta_{2}(\overline{c}))}>0
\]
and
\[
f_{0}(\overline{z},\overline{c})=\frac{\overline{z}\left(  e^{\theta
_{1}(\overline{c})\overline{z}}-e^{\theta_{2}(\overline{c})\overline{z}
}\right)  ^{2}q(\theta_{1}(\overline{c})-\theta_{2}(\overline{c}))}%
{2\theta_{1}(\overline{c}))}>0.
\]

We are going to show that $f_{2}(u,t)<0$ and also
\[
d_{0}(u,t):=g_{1}(u,t)f_{2}(u,t)-f_{1}(u,t)g_{2}(u,t)>0.~
\]
From these inequalities, we conclude that
\[%
\begin{array}
[c]{lll}%
b_{0}(\overline{z},\overline{c})\partial_{x}b_{1}(\overline{z},\overline
{c})+\partial_{xc}b_{0}(\overline{z},\overline{c}) & = & \frac{g_{0}
(\overline{z},\overline{c})}{f_{2}(u,t)}\left(  d_{0}(u,t)+\frac{\partial
_{x}b_{0}(\overline{z},\overline{c})}{f_{0}(\overline{z},\overline{c})}
g_{2}(u,t)\right) \\
& = & \frac{g_{0}(\overline{z},\overline{c})}{f_{2}(u,t)}d_{0}(u,t)\\
& < & 0.
\end{array}
\]
Let us see first that
\[
f_{2}(u,t)=-te^{-4ut}\left(  1+\left(  2u+2ut-1\right)  e^{2u+2ut}\right)
<0.
\]
This holds immediately taking $y=2u+2ut$, because $1+\left(  y-1\right)
e^{y}$ $>0$ for $y>0$. Let us see now that $d_{0}(u,t)>0$, we obtain
\[
d_{0}(u,t)=h_{0}(u,t)(P_{0}(u,t)+P_{1}(u,t)e^{2u}+P_{2}(u,t)e^{2u+2ut}%
+P_{3}(u,t)e^{4u+2ut}+P_{4}(u,t)e^{4u+4ut}),
\]
where%
\[
h_{0}(u,t)=2ut(1+t)e^{-8ut}(e^{2u+2ut}-1)>0\text{ for }t,u>0,
\]

\[
P_{0}(u,t)=t^{2}\text{, }P_{1}(u,t)=(1+t)\left(  2-t+2u+2ut\right)  ,
\]%
\[
P_{2}(u,t)=-2-2u-t+4u^{2}t-t^{2}+4ut^{2}+6u^{2}t^{2}+2ut^{3}-2u^{2}t^{4},
\]%
\[
P_{3}(u,t)=-2+2u-t+2ut+2u^{2}t+t^{2}-2ut^{2}+6u^{2}t^{2}-2ut^{3}+6u^{2}%
t^{3}+2u^{2}t^{4}%
\]
and%
\[
P_{4}(u,t)=2-2u+t-6ut+2u^{2}t-4ut^{2}+4u^{2}t^{2}+2u^{2}t^{3}.
\]
Defining iteratively%
\[%
\begin{array}
[c]{ccc}%
d_{1}(u,t):=\dfrac{\partial_{u}d_{0}(u,t)}{2(1+t)e^{2u}}, &  & d_{2}
(u,t):=\dfrac{\partial_{u}^{2}d_{1}(u,t)}{4e^{2ut}},\\
d_{3}(u,t):=\dfrac{\partial_{u}^{3}d_{2}(u,t)}{8(1+t)^{4}e^{2u}}, &  &
d_{4}(u,t):=\dfrac{\partial_{u}^{3}d_{3}(u,t)}{e^{2ut}};
\end{array}
\]
we obtain%
\begin{equation}
d_{1}(0,t)=\text{ }\partial_{u}d_{1}(0,t)=0, \label{Exp de D1}%
\end{equation}%
\begin{equation}
d_{2}(0,t)=\text{ }\partial_{u}d_{2}(0,t)=\text{ }\partial_{u}^{2}
d_{2}(0,t)=0, \label{Exp de D2}%
\end{equation}%
\begin{equation}
d_{3}(0,t)=5t^{2}\text{, }\partial_{u}d_{3}(0,t)=2t^{2}(7+22t)\text{,
}\partial_{u}^{2}d_{3}(0,t)=8t^{2}(1+18t+25t^{2})\text{, } \label{Exp de D3}%
\end{equation}
and
\[
d_{4}(u,t)=16t^{3}\left(  10+4u+42t+29ut+2u^{2}t+43t^{2}+64ut^{2}+10u^{2}%
t^{2}+44ut^{3}+16u^{2}t^{3}+8u^{2}t^{4}\right)  .
\]
Since $d_{4}(u,t)>0$ and the expressions in (\ref{Exp de D3}) are positive, we
have $d_{3}(u,t)>0$. Similarly, by (\ref{Exp de D2}) and (\ref{Exp de D1}) we
get that $d_{2}(u,t)$, $d_{1}(u,t)\ $and finally $d_{0}(u,t)$ are all
positive. \hfill$\blacksquare$\newline

In the following proposition, we show that the conjecture holds for
$S=[\overline{c}-\varepsilon,\overline{c}]$ with $\varepsilon>0$ small enough.

\begin{proposition}
In the case $\overline{c}$ $>q\sigma^{2}/(2\mu)$, there exists $\varepsilon>0
$ such that $W^{\overline{\zeta}}=V$ in $[0,\infty)\times\lbrack\overline
{c}-\varepsilon,\overline{c}]$, where $\overline{\zeta}$ is the unique
solution of (\ref{Ecuacion diferencial de z0}) with boundary condition
(\ref{Condicion de Optimo en cbarra}).
\end{proposition}

\textit{Proof.} Take $\varepsilon>0$ small enough. By Proposition
\ref{Existencia cerca de Cbarra}, $\overline{\zeta}$ is increasing and so, by
Proposition \ref{z creciente}, $W^{\overline{\zeta}}$ is (2,1)-differentiable
in $[0,\infty)\times\lbrack\overline{c}-\varepsilon,\overline{c}]$. Hence, by
Proposition \ref{Verificacion Curva}, we need to prove that $\partial
_{x}W^{\overline{\zeta}}(\overline{\zeta}(c),c)\leq1$ for $c\in\lbrack
\overline{c}-\varepsilon,\overline{c}]$ and $\partial_{c}W^{\overline{\zeta}%
}(x,c)\leq0$ for $(x,c)\ $with $c\in\lbrack\overline{c}-\varepsilon
,\overline{c}]$ and $0\leq x\leq\overline{\zeta}(c).$

In order to show that $\partial_{x}W^{\overline{\zeta}}(\overline{\zeta
}(c),c)\leq1$ for $c\in\lbrack\overline{c}-\varepsilon,\overline{c}]$, we will
see that $\partial_{x}W^{\overline{\zeta}}(\overline{z},\overline{c})<1$ and
the result will follow by continuity. For the unique point%
\[
x_{0}:=\frac{1}{\theta_{2}(\overline{c})}\log(\frac{-q}{\overline{c}\theta
_{2}(\overline{c})})
\]
where $\partial_{x}W^{\overline{\zeta}}(x_{0},\overline{c})=1$, we have that
$x_{0}>0$ because $\overline{c}$ $>q\sigma^{2}/(2\mu)$.

From Lemma \ref{Punto inicial unico}, the condition $\partial_{x}b_{0}%
(x_{0},\overline{c})<0$ implies $\partial_{x}W^{\overline{\zeta}}(\overline
{z},\overline{c})<1$. Taking $t=\dfrac{-q~}{\theta_{2}(\overline{c}
)\overline{c}}>0$ and $r=-\dfrac{~\theta_{1}(\overline{c})}{\theta
_{2}(\overline{c})}>0$, we can write%
\[
\partial_{x}b_{0}(x_{0},\overline{c})=f(x_{0},\overline{c})g(t,r),
\]
where
\[
f(x_{0},\overline{c})=q^{3}\overline{c}~\theta_{2}^{\prime}(\overline
{c})\theta_{2}^{2}(\overline{c})\left(  e^{\theta_{1}(\overline{c})x_{0}%
}-e^{\theta_{2}(\overline{c})x_{0}}\right)  ^{2}>0,
\]
and%
\[
g(t,r)=-1+\frac{1+(r+1)\log(t)}{t^{r+1}}+\frac{r+1}{r}\left(  1-\frac
{1+(t+1)(r+1)}{t^{r+1}}\right)  .
\]
Since $\overline{c}$ $>q\sigma^{2}/(2\mu)$, we have that $t<1.$ We also get
$g(t,r)<0$, because $g(1,r)=0$ and%
\[
\partial_{t}g(t,r)=\frac{(r+1)^{2}}{t^{r}}(t+1-\log(t))<0.
\]
So $\partial_{x}b_{0}(x_{0},\overline{c})<0.$

Let us show now that $\partial_{c}W^{\overline{\zeta}}(x,c)\leq0$ for
$c\in\lbrack\overline{c}-\varepsilon,\overline{c}]$ and $0\leq x\leq
\overline{\zeta}(c).$ Since%

\[
\partial_{c}W^{\overline{\zeta}}(x,c)=\partial_{c}H^{\overline{\zeta}%
}(x,c)=(e^{\theta_{1}(c)x}-e^{\theta_{2}(c)x})\left(  -b_{0}(x,c)-b_{1}%
(x,c)A^{\overline{\zeta}}(c)+\left(  A^{\overline{\zeta}}\right)  ^{\prime
}(c)\right)  ,
\]
we should analyze the sign of
\[
B(x,c):=-b_{0}(x,c)-b_{1}(x,c)A^{\overline{\zeta}}(c)+\left(  A^{\overline
{\zeta}}\right)  ^{\prime}(c)
\]
for $c\in\lbrack\overline{c}-\varepsilon,\overline{c}]$ and $0\leq
x\leq\overline{\zeta}(c).$ We have that $B(\overline{\zeta}(c),c)=\partial
_{x}B(\overline{\zeta}(c),c)=0.$ Also, from Lemma \ref{Punto inicial unico},
$\partial_{xx}B(\overline{z},\overline{c})=-\partial_{xx}b_{0}(\overline
{z},\overline{c})>0.$ So, $\partial_{xx}B(x,c)>0$ in some neighborhood
\[
U=(\overline{z}-\varepsilon_{1},\overline{z}+\varepsilon_{1})\times
(\overline{c}-\varepsilon_{1},\overline{c}]\subset\lbrack0,\infty
)\times\lbrack\overline{c}-\varepsilon,\overline{c}]
\]
of $(\overline{z},\overline{c})$; and this implies that for any $c\in
\lbrack\overline{c}-\varepsilon,\overline{c}],$ the function $B(\cdot
,c)\ $reaches a strict local maximum at $x=\overline{\zeta}(c)$. In
particular, by Lemma \ref{Punto inicial unico}, $B(\cdot,\overline{c}%
)=-b_{0}(\cdot,\overline{c})+\left(  A^{\overline{\zeta}}\right)  ^{\prime
}(\overline{c})$ reaches the strict global maximum at $x=\overline{z}$ because
$\partial_{x}B(\cdot,\overline{c})=-\partial_{x}b_{0}(\cdot,\overline{c})$
changes from positive to negative at this point. This implies that there
exists a $\delta>0$ such that $B(x,\overline{c})<-\delta$ for $0\leq
x\leq\overline{z}-\varepsilon_{1}.$ Therefore, by continuity arguments, we get
$B(x,c)<0$ for $(x,c)\in\lbrack0,\overline{z}-\varepsilon_{1}]\times
\lbrack\overline{c}-\varepsilon,\overline{c}]$ for some $\varepsilon>0$ small
enough and so we conclude the result. \hfill$\blacksquare$\newline

In the next proposition we show that the optimal value function $V$ is a
uniform limit of $\zeta$-value functions. Moreover, it is a limit of value
functions of extended threshold strategies. The proof uses the convergence
result obtained in Section \ref{Seccion Convergencia}.

\begin{proposition}
Consider, as in Section \ref{Seccion Convergencia}, a sequence of sets
$\mathcal{S}^{n}$ (with $k_{n}$ elements) of the form{\small \ }
\[
\mathcal{S}^{n}=\left\{  c_{1}^{n}=\underline{c}<c_{2}^{n}<\cdots<c_{k_{n}
}^{n}=\overline{c}\right\}
\]
satisfying $\mathcal{S}^{0}=\left\{  \underline{c},\overline{c}\right\}  $,
$\mathcal{S}^{n}\subset\mathcal{S}^{n+1}$ and mesh-size $\delta(\mathcal{S}
^{n}):=\max_{i=2,k_{n}}\left(  c_{i}^{n}-c_{i-1}^{n}\right)  \searrow0$ as $n$
goes to infinity, and the optimal threshold functions $z_{n}^{\ast}
:\widetilde{S}_{n}\rightarrow\lbrack0,\infty)$ defined in Section
\ref{Seccion optima Discreta}. Then, taking $\zeta_{n}(c):=
{\displaystyle\sum\nolimits_{i=1}^{k_{n}-1}}
z_{n}^{\ast}(c_{i}^{n})I_{[c_{i}^{n},c_{i+1}^{n})}$, the $\zeta_{n}$-value
functions $W^{\zeta_{n}}$ converge uniformly to the optimal value function
$V.$
\end{proposition}

\textit{Proof.} Take the functions $V^{n}:[0,\infty)\times\lbrack\underline
{c},\overline{c}]\rightarrow{\mathbb{R}}$ defined in (\ref{Definicion Vn}). By
Proposition \ref{Proposicion Convergencia Uniforme VN a V} and Theorem
\ref{Vbarra es V}, $V^{n}$ converges uniformly to the optimal value function
$V$. Since by Proposition \ref{Proposition Global Lipschitz zone} $V$ is
Lipschitz with constant $K$ and by definition $V^{n}(\cdot,c_{i}^{n}%
)=W^{\zeta_{n}}(\cdot,c_{i}^{n})$ for $c_{i}^{n}\in\mathcal{S}^{n}$, we get
for $c\in\lbrack c_{i}^{n},c_{i+1}^{n})$%

\[%
\begin{array}
[c]{l}%
0\leq V(x,c)-W^{\zeta_{n}}(x,c)\\%
\begin{array}
[c]{ll}%
\leq & \left(  V(x,c)-V(x,c_{i}^{n})\right)  +\left(  V(x,c_{i}^{n}
)-V^{n}(x,c_{i}^{n})\right)  +\left\vert W^{\zeta_{n}}(x,c_{i}^{n}
)-W^{\zeta_{n}}(x,c)\right\vert \\
\leq & K\delta(\mathcal{S}^{n})+\max\left\vert V-V^{n}\right\vert +\left\vert
W^{\zeta_{n}}(x,c_{i}^{n})-W^{\zeta_{n}}(x,c)\right\vert .
\end{array}
\end{array}
\]
Hence, in order to prove the result it suffices to show that there exists a
$K_{1}>0$ such that%
\begin{equation}
\left\vert W^{\zeta_{n}}(x,c_{i}^{n})-W^{\zeta_{n}}(x,c)\right\vert \leq
K_{1}\left\vert c_{i}^{n}-c\right\vert \leq K_{1}\delta(\mathcal{S} ^{n}).
\label{Desigualdad con W}%
\end{equation}
We have that $W^{\zeta_{n}}(x,c_{i+1}^{n})=W^{\zeta_{n}}(x,c)$ for $x\geq
\zeta_{n}(c_{i}^{n})$. So if $\zeta_{n}(c_{i}^{n})>0$ it remains to prove
(\ref{Desigualdad con W}) for $0<x<\zeta_{n}(c_{i}^{n})$. Let us define%
\[
h(x,c):=\dfrac{e^{\theta_{1}(c)x}-e^{\theta_{2}(c)x}}{e^{\theta_{1}
(c)\zeta_{n}(c_{i})}-e^{\theta_{2}(c)\zeta_{n}(c_{i})}}%
\]
and
\[
u(x,c):=\left(  H^{\zeta_{n}}(\zeta_{n}(c_{i}),c_{i+1})-\frac{c}{q}\left(
1-e^{\theta_{2}(c)\zeta_{n}(c_{i})}\right)  \right)  h(x,c).
\]
We can write%
\[%
\begin{array}
[c]{lll}%
\left\vert W^{\zeta_{n}}(x,c_{i}^{n})-W^{\zeta_{n}}(x,c)\right\vert  & = &
\left\vert H^{\zeta_{n}}(x,c_{i}^{n})-H^{\zeta_{n}}(x,c)\ \right\vert \\
& \leq & \left\vert \frac{c_{i}^{n}}{q}\left(  1-e^{\theta_{2}(c_{i}^{n}
)x}\right)  -\frac{c}{q}\left(  1-e^{\theta_{2}(c)x}\right)  \right\vert
+\left\vert u(x,c)-u(x,c_{i})\right\vert .
\end{array}
\]
It is straightforward to see that there exists $K_{1}^{1\text{ }}$such that
\[
\left\vert \frac{c_{i}^{n}}{q}\left(  1-e^{\theta_{2}(c_{i}^{n})x}\right)
-\frac{c}{q}\left(  1-e^{\theta_{2}(c)x}\right)  \right\vert \leq K_{1}%
^{1}\left\vert c-c_{i}^{n}\right\vert \text{ for some }K_{1}^{1}>0.
\]
Since $h(\cdot,c)$ is increasing, we obtain $h(c,0)=0<h(x,c)\leq h(\zeta
_{n}(c_{i}),c)=1$. Also, we have that $0\leq H^{\zeta_{n}}(x,c)\leq
V(x,c)\leq\overline{c}/q\ $for $0<x<\zeta_{n}(c_{i}^{n});$ so using that
$\theta_{2}(c)<0$ it is easy to show that there exist constants $K_{1}%
^{2},K_{1}^{3}>0$ such that
\[%
\begin{array}
[c]{lll}%
\left\vert \partial_{c}u(x,c)\right\vert  & \leq & \left\vert -\frac{1}
{q}\left(  1-e^{\theta_{2}(c)\zeta_{n}(c_{i})}\right)  +\frac{c}{q}\theta
_{2}(c)\zeta_{n}(c_{i})e^{\theta_{2}(c)\zeta_{n}(c_{i})}\frac{\theta
_{2}^{\prime}(c)}{\theta_{2}(c)}\right\vert h(x,c)\\
&  & +\left\vert H^{\zeta_{n}}(\zeta_{n}(c_{i}),c_{i+1})-\frac{c}{q}\left(
1-e^{\theta_{2}(c)\zeta_{n}(c_{i})}\right)  \right\vert \left\vert
\partial_{c}h(x,c)\right\vert \\
& \leq & K_{1}^{2}+K_{1}^{3}\left\vert \partial_{c}h(x,c)\right\vert .
\end{array}
\]
Calling $y=\zeta_{n}(c_{i})>0\ $and $\rho=\frac{x}{\zeta_{n}(c_{i})}\in(0,1),
$ we obtain%
\[%
\begin{array}
[c]{lll}%
\partial_{c}h(x,c) & = & T(\rho,y,c)\\
& = & \dfrac{\left(  \theta_{1}^{\prime}(c)\rho ye^{\theta_{1}(c)\rho
y}-\theta_{2}^{\prime}(c)\rho ye^{\theta_{2}(c)\rho y}\right)  \left(
e^{\theta_{1}(c)y}-e^{\theta_{2}(c)y}\right)  }{\left(  e^{\theta_{1}
(c)y}-e^{\theta_{2}(c)y}\right)  ^{2}}\text{ }\\
&  & -\dfrac{\left(  \theta_{1}^{\prime}(c)ye^{\theta_{1}(c)y}-\theta
_{2}^{\prime}(c)ye^{\theta_{2}(c)y}\right)  \left(  e^{\theta_{1}(c)\rho
y}-e^{\theta_{2}(c)\rho y}\right)  }{\left(  e^{\theta_{1}(c)y}-e^{\theta
_{2}(c)y}\right)  ^{2}}\text{ .}%
\end{array}
\]
Now, on the one hand,
\[
T(0,y,c)=T(1,y,c)=0\text{, }\lim_{y\rightarrow0}T(\rho,y,c)=0
\]
and on the other hand, taking $\varepsilon>0$, and $y\geq\varepsilon,$ there
exists $K_{1}^{4}>0$ such that%
\[%
\begin{array}
[c]{lll}%
T(\rho,y,c) & = & -y(1-\rho)\theta_{1}(c)e^{-y(1-\rho)\theta_{1}(c)}
\dfrac{\theta_{1}^{\prime}(c)}{\left(  1-e^{(\theta_{2}(c)-\theta_{1}
(c))y}\right)  ^{2}\theta_{1}(c)}\\
&  & +\theta_{1}(c)ye^{-\theta_{1}(c)y}\dfrac{e^{\theta_{2}(c)\rho y}-\rho
e^{\left(  (\rho-1)\theta_{1}(c)+\theta_{2}(c)\right)  y}}{\left(
1-e^{(\theta_{2}(c)-\theta_{1}(c))y}\right)  ^{2}}\frac{\theta_{1}^{\prime
}(c)}{\theta_{1}(c)}\\
&  & +\theta_{1}(c)ye^{-\theta_{1}(c)y}\dfrac{e^{\left(  (\rho-1)\theta
_{1}(c)+\theta_{2}(c)\right)  y}+e^{\left(  (\rho+1)\theta_{2}(c)-\theta
_{1}(c)\right)  y}(\rho-1)-\rho e^{\theta_{2}(c)\rho y}}{\left(
1-e^{(\theta_{2}(c)-\theta_{1}(c))y}\right)  ^{2}}\frac{\theta_{2}^{\prime
}(c)}{\theta_{1}(c)}\\
& \leq & K_{1}^{4},
\end{array}
\]
because $se^{-s}$ is bounded for $s\geq0.$ So we get (\ref{Desigualdad con W})
and finally the result.\hfill$~\blacksquare$

\section{Numerical examples \label{Numerical examples}}

Let us finally consider a numerical illustration for the case $\mu=4$,
$\sigma=2$ and $q=0.1$ for $S=[0,\overline{c}]$. In order to obtain the
corresponding optimal value function $V^{S},$ we proceed as follows:

\begin{enumerate}
\item We obtain $\overline{\zeta}$ solving numerically
the ordinary differential equation (\ref{Ecuacion diferencial de z0}) with
boundary condition (\ref{Condicion de Optimo en cbarra}), using the Euler method.

\item We check that the $\overline{\zeta}$-value function $W^{\overline{\zeta
}}$ defined in (\ref{Definicion de Wz}) satisfies the conditions of
Proposition \ref{Verificacion Curva}. Hence $W^{\overline{\zeta}}$ is the
optimal value function $V^{S}$.
\end{enumerate}

\begin{figure}
	\centering
	\begin{subfigure}[b]{0.4\textwidth}
		\centering
		\includegraphics[width=2.5in]{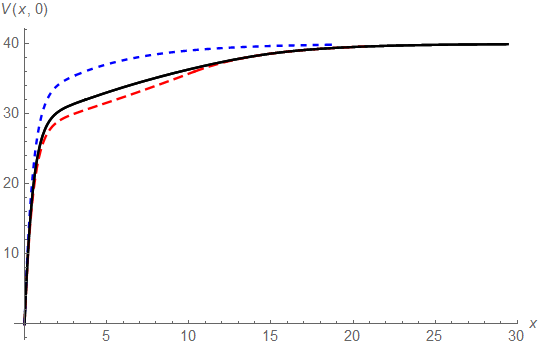}
		\caption{{\small $V^S(x,0)$ (black) together with $V_{NR}(x)$ (blue) and $V_{1}(x)$ (red)}}
		\label{fig1}
	\end{subfigure}
	\hspace*{0.2cm}
	\begin{subfigure}[b]{0.4\textwidth}
		\centering
		\includegraphics[width=2.5in]{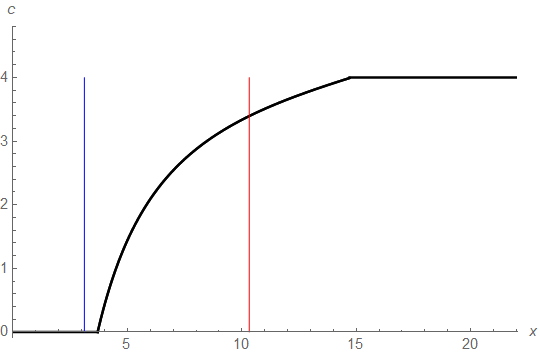}
		\caption{{\small Optimal curve $(\zeta(c),c)$ (black) together with $b^*$ (blue) and $b_R^*$ (red) }}
		\label{fig2}
	\end{subfigure}\caption{$\overline{c}=4$}
	\end{figure}

\begin{figure}
	\centering
	\begin{subfigure}[b]{0.4\textwidth}
		\centering
		\includegraphics[width=2.5in]{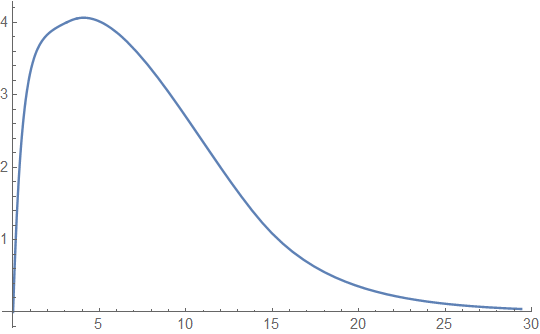}
		\caption{{\small $V_{NR}(x)-V^S(x,0)$ as a function of $x$ ($\overline{c}=4$)}}
		\label{fig3}
	\end{subfigure}
\hspace*{0.2cm}
	\begin{subfigure}[b]{0.4\textwidth}
		\centering
		\includegraphics[width=2.5in]{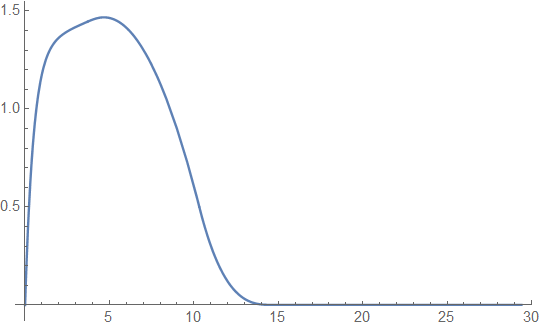}
		\caption{{\small $V^S(x,0)-V_1(x)$ as a function of $x$ ($\overline{c}=4$)}}
		\label{fig4}
	\end{subfigure}\caption{$\overline{c}=4$}
\end{figure}

Let us first consider the case $\overline{c}=4$ (i.e.\ the maximal allowed
dividend rate is the drift of the surplus process $X_{t}$). Figure \ref{fig1}
depicts $V^{S}(x,0)$ as a function of initial capital $x$ together with the
value function $V_{NR}(x)$ of the classical dividend problem without
ratcheting constraint, for which the optimal strategy is a threshold strategy
of not paying any dividends when the surplus level is below $b^{*}$ and pay
dividends at rate $\overline{c}$ above $b^{*}$. Recall from Asmussen and
Taksar \cite{asmtak} or also Gerber and Shiu \cite{Gerber2006} that in the
notation of the present paper
\[
V_{NR}(x)=\left\{
\begin{array}
[c]{ll}%
\frac{\overline{c}}{q}\,\frac{e^{\theta_{1}(0)x}-e^{\theta_{2}(0)x}}%
{\theta_{1}(0)\,e^{\theta_{1}(0)b^{*}}-\theta_{2}(0)\,e^{\theta_{2}(0)b^{*}}%
}, & 0\le x\le b^{*},\\
\frac{\overline{c}}{q}+e^{\theta_{2}(\overline{c})(x-b^{*})}/\theta
_{2}(\overline{c}), & x\ge b^{*}%
\end{array}
\right.
\]
with optimal threshold
\[
b^{*}=\frac{1}{\theta_{1}(0)-\theta_{2}(0)}\log\frac{\theta_{2}(0)\;(\theta
_{2}(0)-\theta_{2}(\overline{c}))}{\theta_{1}(0)\;(\theta_{1}(0)-\theta
_{2}(\overline{c}))}.
\]
One observes that for both small and large initial capital $x$ the efficiency
loss when introducing the ratcheting constraint is very small, only for
intermediate values of $x$ the resulting expected discounted dividends are
significantly smaller, but even there the relative efficiency loss is not big
(see Figure \ref{fig3} for a plot of this difference). We also compare
$V^{S}(x,0)$ in Figure \ref{fig1} with the optimal value function
\[
V_{1}(x):=V^{0}(x)\;\text{for }S=\{0,\overline{c}\}
\]
of the further constrained one-step ratcheting strategy, where only once
during the lifetime of the process the dividend rate can be increased from 0
to $\overline{c}$. That latter case was studied in detail in \cite{ABB}, where
it was also shown that the optimal threshold level $b_{R}^{*}$ for that switch
is exactly the one for which the resulting expected discounted dividends match
with the ones of a threshold strategy underlying $V_{NR}$, but at the (for the
latter problem non-optimal) threshold $b_{R}^{*}$. We observe that the
performance of this simple one-step ratcheting is already remarkably close to
the one of the overall optimal ratcheting strategy represented by $V^{S}(x,0)$
(see also Figure \ref{fig4} for a plot of the difference). A similar effect
had already been observed for the optimal ratcheting in the Cram\'er-Lundberg
model (cf.\ \cite{AAM17}).

Figure \ref{fig2} depicts the optimal ratcheting curve $(\zeta(c),c)$
underlying $V^{S}(x,0)$ for this example together with the optimal threshold
$b^{\ast}$ of the unconstrained dividend problem and the optimal switching
barrier $b_{R}^{\ast}$ for the one-step ratcheting strategy. One sees that the
irreversibility of the dividend rate increase in the ratcheting case leads to
a rather conservative behavior of not starting any (even not small) dividend
payments until the surplus level is above the optimal threshold level
$b^{\ast}$ underlying the non-constrained dividend problem. On the other hand,
the one-step ratcheting strategy with optimal switching barrier $b_{R}^{\ast}$
roughly in the middle of the optimal curve already leads to a remarkably good
approximation (lower bound) for the performance of the overall optimal
ratcheting strategy. \newline

\begin{figure}
	\centering
	\begin{subfigure}[b]{0.4\textwidth}
		\centering
		\includegraphics[width=2.5in]{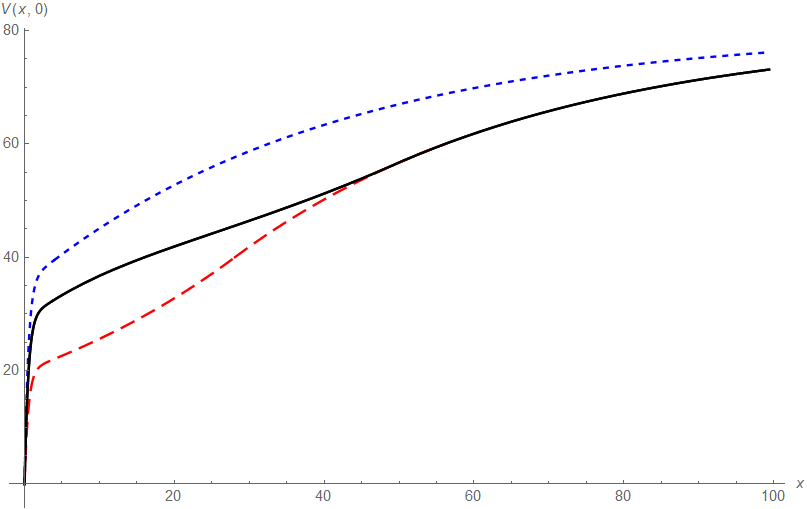}
		\caption{{\small $V^S(x,0)$ (black) together with $V_{NR}(x)$ (blue) and $V_1(x)$ (red)}}
		\label{fig5}
	\end{subfigure}\hspace*{0.2cm}
	\begin{subfigure}[b]{0.4\textwidth}
		\centering
		\includegraphics[width=2.5in]{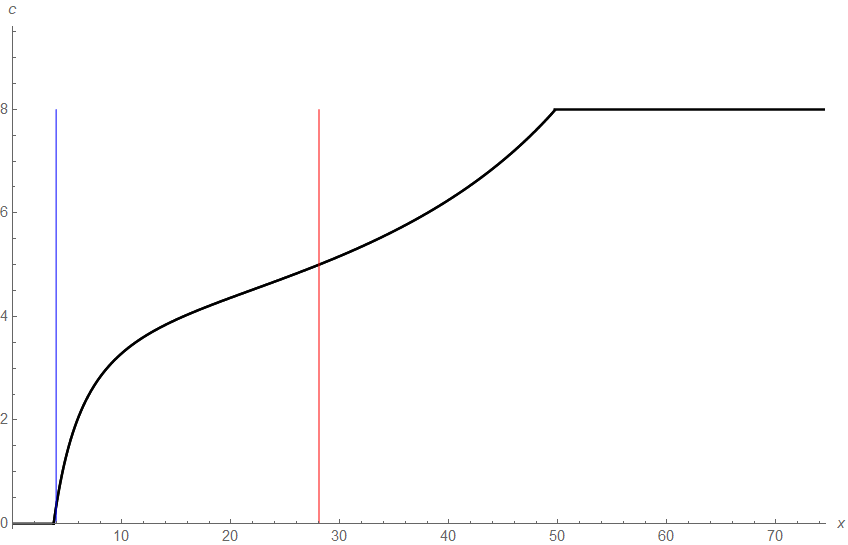}
		\caption{{\small Optimal curve $(\zeta(c),c)$ (black) together with $b^*$ (blue) and $b_R^*$ (red)}}
		\label{fig6}
	\end{subfigure}\caption{$\overline{c}=8$}\label{figg3}
\end{figure}

\begin{figure}
	\centering
	\begin{subfigure}[b]{0.4\textwidth}
		\centering
		\includegraphics[width=2.5in]{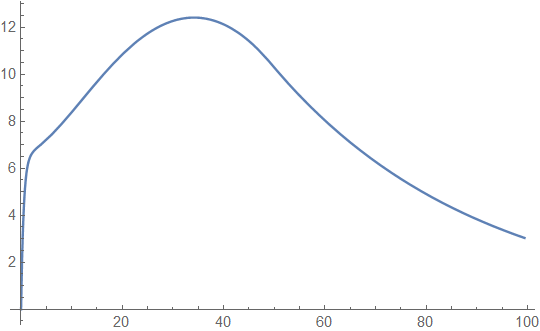}
		\caption{{\small $V^S(x,0)$ (black) together with $V_{NR}(x)$ (blue) and $V_1(x)$ (red)}}
		\label{fig7}
	\end{subfigure}
\hspace*{0.2cm}
	\begin{subfigure}[b]{0.4\textwidth}
		\centering
		\includegraphics[width=2.5in]{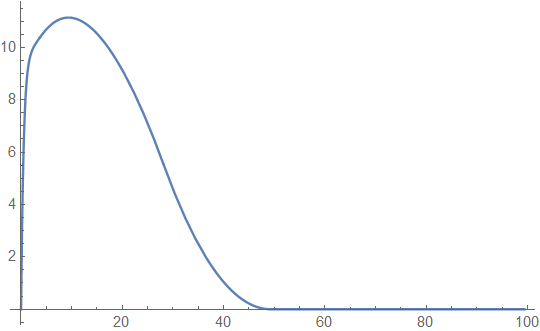}
		\caption{{\small $V^S(x,0)-V_1(x)$ as a function of $x$ }}
		\label{fig8}
	\end{subfigure}\caption{$\overline{c}=8$}\label{figg4}
\end{figure}

In Figures \ref{figg3} and \ref{figg4} we give the analogous plots for the case
$\overline{c}=8$, so that the maximal dividend rate is twice as large as the
drift $\mu$ of the uncontrolled risk process. The global picture is quite
similar, also in this case the efficiency loss introduced by ratcheting is
more pronounced and present also for larger initial capital $x$. Also, the
further efficiency loss by restricting to a simple one-step ratcheting
strategy is considerably larger for not too large $x$. Finally, in that case
the first increase of dividends already happens for surplus values (slightly)
smaller than the optimal threshold $b^{*}$ of the unconstrained case.

\section{Conclusion}

\label{seccon} In this paper we studied and solved the problem of finding
optimal dividend strategies in a Brownian risk model, when the dividend rate
can not be decreased over time. We showed that the value function is the
unique viscosity solution of a two-dimensional Hamilton-Jacobi-Bellman
equation and it can be approximated arbitrarily closely by threshold
strategies for finitely many possible dividend rates, which are established to
be optimal in their discrete setting. We used calculus of variation techniques
to identify the optimal curve that separates the state space into a change and
a non-change region and provided partial results for the overall optimality of
this strategy (which can be seen as a two-dimensional analogue of the
optimality of dividend threshold strategies in the one-dimensional diffusion
setting without the ratcheting constraint). In contrast to \cite{AAM}, the
same analysis is applicable for all finite levels of maximal dividend rate
$\overline{c}$, i.e.\ also if the latter exceeds the drift $\mu$. We also gave
some numerical examples determining the optimal curve strategy. These results
illustrate that the ratcheting constraint does not reduce the efficiency of
the optimal dividend strategy substantially and that, much as in the compound
Poisson setting, the simpler strategy of only stepping up the dividend rate
once during the lifetime of the process is surprisingly close to optimal in
absolute terms. In terms of a possible direction of future research, as
mentioned in Section \ref{Optimal strategies} we conjecture that a curve
strategy dividing the state space into a change and a non-change region is
optimal in full generality for the diffusion model, and it remains open to
formally prove the latter. Furthermore, it could be interesting to extend the
results of the present paper to the case where the dividend rate may be
decreased by a certain percentage of its current value (see e.g.\ \cite{bayr})
or to place the dividend consumption pattern into a general habit formation
framework (see e.g.\ \cite{bayr2} for an interesting related paper in a
deterministic setup).

\section{Appendix}

\label{secapp} \textit{Proof of Proposition
\ref{Proposition Global Lipschitz zone}.} By Proposition
\ref{Monotone Optimal Value Function}, we have
\begin{equation}
0\leq V^{S}(x_{2},c_{1})-V^{S}(x_{1},c_{2}) \label{Lips1}%
\end{equation}
for all $0\leq x_{1}\leq x_{2}$ and $c_{1},c_{2}\in S$ with $c_{1}\leq c_{2}.
$

Let us show now, that there exists $K_{1}>0$ such that%
\begin{equation}
V^{S}(x_{2},c)-V^{S}(x_{1},c)\leq K_{1}\left(  x_{2}-x_{1}\right)
\label{Lips2}%
\end{equation}
for all $0\leq x_{1}\leq x_{2}$. Take $\varepsilon>0$ and $C\in\Pi_{x_{2}
,c}^{S}$ such that%

\begin{equation}
J(x_{2};C)\geq V^{S}(x_{2},c)-\varepsilon, \label{casi optima x2}%
\end{equation}
the associated control process is given by%

\[
X_{t}^{C}=x_{2}+\int_{0}^{t}(\mu-C_{s})ds+W_{t}.
\]
Let $\tau$ be the ruin time of the process $X_{t}^{C}$. Define $\widetilde
{C}\in\Pi_{x_{1},c}^{S}$ as $\widetilde{C}_{t}=C_{t}$ and the associated
control process
\[
\text{ }X_{t}^{\widetilde{C}}=x_{1}+\int_{0}^{t}(\mu-C_{s})ds+W_{t}.
\]
Let $\widetilde{\tau}\leq\tau$ be the ruin time of the process $X_{t}%
^{\widetilde{C}}$; it holds that $X_{t}^{C}-X_{t}^{\widetilde{C}}=x_{2}-x_{1}$
for $t\leq\widetilde{\tau}$. Hence we have
\begin{equation}%
\begin{array}
[c]{lll}%
V^{S}(x_{2},c)-V^{S}(x_{1},c) & \leq & J(x_{2};C)-J(x_{1};\widetilde
{C})+\varepsilon\\
& \leq & V^{S}(x_{2}-x_{1},0)+\varepsilon\\
& \leq & V_{NR}(x_{2}-x_{1})+\varepsilon\\
& \leq & K_{1}(x_{2}-x_{1})+\varepsilon.
\end{array}
\label{diferencia de vs}%
\end{equation}
So, by Remark \ref{Optima sin ratcheting}, we have (\ref{Lips2}) with
$K_{1}=V_{NR}^{\prime}(0).$

Let us show now that, given $c_{1},c_{2}\in S$ with $c_{1}\leq c_{2},$ there
exists $K_{2}>0$ such that%
\begin{equation}
V^{S}(x,c_{1})-V^{S}(x,c_{2})\leq K_{2}\left(  c_{2}-c_{1}\right)  .
\label{Lips3}%
\end{equation}
Take $\varepsilon>0$ and $C\in\Pi_{x,c_{1}}^{S}$ such that%

\begin{equation}
J(x;C)\geq V^{S}(x,c_{1})-\varepsilon, \label{casiOptimac1}%
\end{equation}
define the stopping time
\begin{equation}
\widehat{T}=\min\{t:C_{t}\geq c_{2}\} \label{Definicion TauSombrero}%
\end{equation}
and denote $\tau$ the ruin time of the process $X_{t}^{C}$. Let us consider
$\widetilde{C}\in\Pi_{x,c_{2}}^{S}\ $as $\widetilde{C}_{t}=c_{2}%
I_{t<\widehat{T}}+C_{t}I_{t\geq\widehat{T}}$; denote by $X_{t}^{\widetilde{C}%
}$ the associated controlled surplus process and by $\overline{\tau}\leq\tau$
the corresponding ruin time. We have that $\widetilde{C}_{s}-C_{s}\leq
c_{2}-c_{1}$ and so $X_{\overline{\tau}}^{C}=X_{\overline{\tau}}%
^{C}-X_{\overline{\tau}}^{\widetilde{C}}\leq(c_{2}-c_{1})\overline{\tau}$,
which implies%
\[
\int_{\overline{\tau}}^{\tau}C_{s}e^{-q\left(  s-\overline{\tau}\right)
}ds\leq V_{NR}((c_{2}-c_{1})\overline{\tau}).
\]
Hence, we can write,%
\begin{equation}%
\begin{array}
[c]{lll}%
V^{S}(x,c_{1})-V^{S}(x,c_{2}) & \leq & J(x;C)+\varepsilon-J(x;\widetilde{C})\\
& = & \mathbb{E}\left[  \int_{0}^{\overline{\tau}}\left(  C_{s}-\widetilde
{C}_{s}\right)  e^{-qs}ds\right]  +\mathbb{E}\left[  \int_{\overline{\tau}
}^{\tau}C_{s}e^{-qs}ds\right]  +\varepsilon\\
& \leq & 0+\mathbb{E}\left[  \int_{\overline{\tau}}^{\tau}C_{s}e^{-qs}
ds\right]  +\varepsilon\\
& \leq & E[e^{-q\overline{\tau}}\int_{\overline{\tau}}^{\tau}C_{s}e^{-q\left(
s-\overline{\tau}\right)  }ds]+\varepsilon\\
& \leq & K_{1}E[e^{-q\overline{\tau}}\overline{\tau}(c_{2}-c_{1}
)]+\varepsilon\\
& \leq & K_{2}(c_{2}-c_{1})+\varepsilon.
\end{array}
\label{Lipaux1C}%
\end{equation}
So, we deduce (\ref{Lips3}), taking $K_{2}=K_{1}\max_{t\geq0}\{e^{-qt}t\}.$ We
conclude the result from (\ref{Lips1}), (\ref{Lips2}) and (\ref{Lips3}).
$\blacksquare$\newline

\textit{Proof of Proposition \ref{Proposicion Viscosidad}.} Let us show first
that $V$ is a viscosity supersolution in $(0,\infty)\times\lbrack\underline
{c},\overline{c})$ . By Proposition \ref{Monotone Optimal Value Function},
$\partial_{c}V\leq0$ in $(0,\infty)\times\lbrack\underline{c},\overline{c})$
in the viscosity sense.

Consider now $(x,c)\in(0,\infty)\times\lbrack\underline{c},\overline{c})$ and
the admissible strategy $C\in\Pi_{x,c}^{S}$, which pays dividends at constant
rate $c$ up to the ruin time $\tau$. Let $X_{t}^{C}$ be the corresponding
controlled surplus process and suppose that there exists a test function
$\varphi$ for supersolution (\ref{HJB equation}) at $(x,c)$. Using Lemma
\ref{DPP}, we get for $h>0$%

\[%
\begin{array}
[c]{lll}%
\varphi(x,c) & = & V(x,c)\\
& \geq & \mathbb{E}\left[  \int\nolimits_{0}^{\tau\wedge h}e^{-q\,s}
\,cds\right]  +\mathbb{E}\left[  e^{-q(\tau\wedge h)}\varphi(X_{\tau\wedge
h}^{C},c))\right]  .
\end{array}
\]
Hence, using It\^{o}'s formula%
\[%
\begin{array}
[c]{lll}%
0 & \geq & \mathbb{E}\left[  \int\nolimits_{0}^{\tau\wedge h}e^{-q\,s}
\,c\,ds\right]  +\mathbb{E}\left[  I_{\tau>h}\left(  e^{-q\,h}\varphi
(X_{s}^{C},c)-\varphi(x,c)\right)  \right]  -\varphi(x,c)\mathbb{P(}h>\tau)\\
& = & \mathbb{E}\left[  \int\nolimits_{0}^{\tau\wedge h}e^{-q\,s}
\,c\,ds\right]  +\mathbb{E}\left[  I_{\tau>h}\int\nolimits_{0}^{h}
e^{-q\,s}(\frac{\sigma^{2}}{2}\partial_{xx}\varphi(X_{s}^{C},c)+\partial
_{x}\varphi(X_{s}^{C},c)(\mu-c)-q\varphi(X_{s}^{C},c))ds\right]
-\varphi(x,c)\mathbb{P(}h>\tau).
\end{array}
\]
So, dividing by $h$ and taking $h\rightarrow0^{+}$, we get $\mathcal{L}%
^{c}(\varphi)(x,c)\leq0;$, so that $V$ is a viscosity supersolution at $(x,c)$.

Let us prove now that $V$ it is a viscosity subsolution in $(0,\infty
)\times\lbrack\underline{c},\overline{c})$. Assume first that $V$ is not a
subsolution of (\ref{HJB equation}) at $\left(  x,c\right)  \in(0,\infty
)\times\lbrack\underline{c},\overline{c})$. Then there exist $\varepsilon>0$,
$0<h<\min\left\{  x/2,\overline{c}-c\right\}  $ and a (2,1)-differentiable
function $\psi$ with $\psi(x,c)=V(x,c)$ such that $\psi\geq V$,%

\begin{equation}
\max\{\mathcal{L}^{c}(\psi)(y,d),\partial_{c}\psi(y,d)\}\leq-q\varepsilon<0
\label{Desig1aprima}%
\end{equation}
for $\left(  y,d\right)  \in$ $[x-h,x+h]\times\lbrack c,c+h]$ and%

\begin{equation}
V(y,d)\leq\psi(y,d)-\varepsilon\label{Desig2a}%
\end{equation}
for $\left(  y,d\right)  \notin\lbrack x-h,x+h]\times\lbrack c,c+h]$. Consider
the controlled risk process $X_{t}$ corresponding to an admissible strategy
$C\in\Pi_{x,c}^{S}$ and define
\[
\tau^{\ast}=\inf\{t>0:\text{ }\left(  X_{t},C_{t}\right)  \notin\lbrack
x-h,x+h]\times\lbrack c,c+h]\}\text{.}%
\]
Since $C_{t}$ is non-decreasing and right-continuous, it can be written as
\begin{equation}
C_{t}=c+\int\nolimits_{0}^{t}dC_{s}^{co}+\sum_{\substack{C_{s}\neq C_{s^{-}}
\\0\leq s\leq t}}(C_{s}-C_{s^{-}}), \label{Lt para ITO}%
\end{equation}
where $C_{s}^{co}$ is a continuous and non-decreasing function.

Take a (2,1)-differentiable function $\psi:(0,\infty)\times\lbrack
\underline{c},\overline{c}]\rightarrow\lbrack0,\infty)$. Using the expression
(\ref{Lt para ITO}) and the change of variables formula (see for instance
\cite{Protter}), we can write
\begin{equation}%
\begin{array}
[c]{l}%
e^{-q\tau^{\ast}}\psi(X_{\tau^{\ast}}^{C},C_{\tau^{\ast}})-\psi(x,c)\\%
\begin{array}
[c]{ll}%
= & \int\nolimits_{0}^{\tau^{\ast}}e^{-qs}\partial_{x}\psi(X_{s}^{C},C_{s^{-}
})(\mu-C_{s^{-}})ds+\int\nolimits_{0}^{\tau^{\ast}}e^{-qs}\partial_{c}
\psi(X_{s}^{C},C_{s^{-}})dC_{s}^{co}\\
& +\sum_{\substack{C_{s}\neq C_{s^{-}} \\0\leq s\leq\tau^{\ast}}}e^{-qs}
(C_{s}-C_{s^{-}})\partial_{c}\psi(X_{s}^{C},C_{s^{-}})\\
& +\int\nolimits_{0}^{\tau^{\ast}}e^{-qs}(-q\psi(X_{s}^{C},C_{s^{-}}
)+\frac{\sigma^{2}}{2}\partial_{xx}\psi(X_{s}^{C},C_{s^{-}}))ds+\int
\nolimits_{0}^{\tau^{\ast}}e^{-qs}\partial_{x}\psi(X_{s}^{C},C_{s^{-}})\sigma
dW_{s}.
\end{array}
\end{array}
\label{Paso 1}%
\end{equation}
Hence, from (\ref{Desig1aprima}), we can write%
\[%
\begin{array}
[c]{l}%
\mathbb{E}\left[  e^{-q\tau^{\ast}}\psi(X_{\tau^{\ast}}^{C},C_{\tau^{\ast}
})\right]  -\psi(x,c)\\%
\begin{array}
[c]{ll}%
= & \mathbb{E}\left[  \int\nolimits_{0}^{\tau^{\ast}}e^{-qs}\mathcal{L}
^{C_{s^{-}}}(\psi)(X_{s^{-}}^{C},C_{s^{-}})ds-\int\nolimits_{0}^{\tau^{\ast}
}e^{-qs}C_{s^{-}}ds\right] \\
& +\mathbb{E}\left[  \int\nolimits_{0}^{\tau^{\ast}}e^{-qs}\partial_{c}
\psi(X_{s^{-}}^{C},C_{s^{-}})dC_{s}^{c}+\sum_{\substack{C_{s}\neq C_{s^{-}}
\\0\leq s\leq\tau^{\ast}}}e^{-qs}(C_{s}-C_{s^{-}})\partial_{c}\psi(X_{s^{-}
}^{C},C_{s^{-}})\right] \\
\leq & \mathbb{E}\left[  \varepsilon\left(  e^{-q\tau^{\ast}}-1\right)
-\int\nolimits_{0}^{\tau^{\ast}}e^{-qs}C_{s^{-}}ds-q\varepsilon\left(
\int\nolimits_{0}^{\tau^{\ast}}e^{-qs}dC_{s}\right)  \right]  .
\end{array}
\end{array}
\]

So, from (\ref{Desig2a})%
\[%
\begin{array}
[c]{l}%
\mathbb{E}\left[  e^{-q\tau^{\ast}}V(X_{\tau^{\ast}}^{C},C_{\tau^{\ast}
})\right] \\%
\begin{array}
[c]{cl}%
\leq & \mathbb{E}\left[  \psi(x,c)-e^{-q\tau^{\ast}}\varepsilon\right]
+\mathbb{E}\left[  \psi(X_{\tau^{\ast}}^{C},C_{\tau^{\ast}})e^{-q\tau^{\ast}
}-\psi(x,c)\right] \\
\leq & \psi(x,c)-\varepsilon-\mathbb{E}(\int\nolimits_{0}^{\tau^{\ast}}
e^{-qs}C_{s^{-}}ds).
\end{array}
\end{array}
\]
Hence, using Lemma \ref{DPP}, we have that
\[
V(x,c)=\sup\limits_{C\in\Pi_{x,c}^{S}}\mathbb{E}\left(  \int\nolimits_{0}%
^{\tau^{\ast}}e^{-qs}C_{s^{-}}ds+e^{-c\tau^{\ast}}V(X_{\tau^{\ast}}%
^{C},C_{\tau^{\ast}})\right)  \leq\psi(x,c)-\varepsilon.
\]
but this is a contradiction because we have assumed that $V(x,c)=\psi(x,c)$.
So we have the result. $\blacksquare$ \bigskip

\textit{Proof of Lemma \ref{Lema para Unicidad}.} A locally Lipschitz function
$\overline{u}$\ $:[0,\infty)\times\lbrack\underline{c},\overline
{c}]\rightarrow{\mathbb{R}}$\ is a viscosity supersolution of
(\ref{HJB equation}) at $(x,c)\in(0,\infty)\times(\underline{c},\overline{c}%
)$, if any test function $\varphi$ for supersolution at $(x,c)$ satisfies
\begin{equation}
\max\{\mathcal{L}^{c}(\varphi)(x,c),\partial_{c}\varphi(x,c)\}\leq0\text{,}
\label{superyalfa}%
\end{equation}
and a locally Lipschitz function $\underline{u}:[0,\infty)\times
\lbrack\underline{c},\overline{c}]\rightarrow{\mathbb{R}}$\ is a viscosity
subsolution\ of (\ref{HJB equation}) at $(x,c)\in(0,\infty)\times
(\underline{c},\overline{c})$\ if any test function $\psi$ for subsolution at
$(x,c)$ satisfies
\begin{equation}
\max\{\mathcal{L}^{c}(\psi)(x,c),\partial_{c}\psi(x,c)\}\geq0.
\label{subxalfa}%
\end{equation}

Suppose that there is a point $(x_{0},c_{0})\in\lbrack0,\infty)\times
(\underline{c},\overline{c})$ such that $\underline{u}(x_{0},c_{0}%
)-\overline{u}(x_{0},c_{0})>0$. Let us define $h(c)=1+e^{-{c}/ {{\overline{c}%
}}}$ and
\[
\overline{u}^{s}(x,c)=s\,h(c)\,\overline{u}(x,c)
\]
for any $s>1$. We have that $\varphi$ is a test function for supersolution of
$\overline{u}$ at $(x,c)$ if and only if $\varphi^{s}=s\,h(c)\,\varphi$ is a
test function for supersolution of $\overline{u}^{s}$ at $(x,c)$. We have
\begin{equation}
\mathcal{L}^{c}(\varphi^{s})(x,c)=sh(c)\mathcal{L}^{c}(\varphi
)(x,c)+c(1-sh(c))<0, \label{Desigualdad L1}%
\end{equation}
and%
\begin{equation}
\partial_{c}\varphi^{s}(x,c)\leq-\frac{s}{\overline{c}}\varphi(x,c)e^{-\frac
{c}{\overline{c}}}<0 \label{Desigualdad L2}%
\end{equation}
for $\varphi(x,c)>0.$ Take $s_{0}>1$ such that $\underline{u}(x_{0}%
,c_{0})-\overline{u}^{s_{0}}(x_{0},c_{0})>0$. We define%
\begin{equation}
M=\sup\limits_{x\geq0,\underline{c}\leq c\leq\overline{c}}\left(
\underline{u}(x,c)-\overline{u}^{s_{0}}(x,c)\right)  . \label{Definicion de M}%
\end{equation}
Since $\lim_{x\rightarrow\infty}\underline{u}(x,c)\leq\overline{c}/q\leq
\lim_{x\rightarrow\infty}\overline{u}(x,c)$, there exist $b>x_{0}$ such that%

\begin{equation}
\sup\limits_{\underline{c}\leq c\leq\overline{c}}\underline{u}(x,c)-\overline
{u}^{s_{0}}(x,c)<0\text{ for }x\geq b. \label{desborde}%
\end{equation}
From (\ref{desborde}), we obtain that%

\begin{equation}
0<\underline{u}(x_{0},c_{0})-\overline{u}^{s_{0}}(x_{0},c_{0})\leq
M:=\max\limits_{x\in\left[  0,b\right]  ,\underline{c}\leq c\leq\overline{c}
}\left(  \underline{u}(x,c)-\overline{u}^{s_{0}}(x,c)\right)  . \label{desmax}%
\end{equation}
Call $\left(  x^{\ast},c^{\ast}\right)  :=\arg\max\limits_{x\in\left[
0,b\right]  ,\underline{c}\leq c\leq\overline{c}}\left(  \underline
{u}(x,c)-\overline{u}^{s_{0}}(x,c)\right)  $. Let us consider the set%

\[
\mathcal{A}=\left\{  \left(  x,y,c,d\right)  :0\leq x\leq y\leq b\text{,
}\underline{c}\leq\ c\leq\overline{c}\text{, }\underline{c}\leq d\leq
\overline{c}\right\}
\]
and, for all $\lambda>0$, the functions%

\begin{equation}%
\begin{array}
[c]{l}%
\Phi^{\lambda}\left(  x,y,c,d\right)  =\dfrac{\lambda}{2}\left(  x-y\right)
^{2}+\dfrac{\lambda}{2}\left(  c-d\right)  ^{2}+\frac{2m}{\lambda^{2}\left(
y-x\right)  +\lambda},\\
\Sigma^{\lambda}\left(  x,y,c,d\right)  =\underline{u}(x,c)-\overline
{u}^{s_{0}}(y,d)-\Phi^{\lambda}\left(  x,y,c,d\right)  .
\end{array}
\label{def-sigma-definitiva-por-hoy}%
\end{equation}
$\allowbreak$ Calling $M^{\lambda}=\max\limits_{A}\Sigma^{\lambda}$ and
$\left(  x_{\lambda},y_{\lambda},c_{\lambda},d_{\lambda}\right)  =\arg
\max\limits_{A}\Sigma^{\lambda}$, we obtain that $M^{\lambda}\geq
\Sigma^{\lambda}(x^{\ast},x^{\ast},c^{\ast},c^{\ast})=M-\frac{2m}{\lambda}$,
and so%

\begin{equation}
\liminf\limits_{\lambda\rightarrow\infty}M^{\lambda}\geq M.
\label{liminf-mlambda}%
\end{equation}

There exists $\lambda_{0}$ large enough such that if $\lambda\geq\lambda_{0}$,
then $\left(  x_{\lambda},y_{\lambda},c_{\lambda},d_{\lambda}\right)  $
$\notin\partial A$, the proof is similar to the one of Lemma 4.5 of
{\cite{AAM}.}

Using the inequality%

\[
\Sigma^{\lambda}\left(  x_{\lambda},x_{\lambda},c_{\lambda},c_{\lambda
}\right)  +\Sigma^{\lambda}\left(  y_{\lambda},y_{\lambda},d_{\lambda
},d_{\lambda}\right)  \leq2\Sigma^{\lambda}\left(  x_{\lambda},y_{\lambda
},c_{\lambda},d_{\lambda}\right)  ,
\]
we obtain that%

\[
\lambda\left\Vert (x_{\lambda}-y_{\lambda},c_{\lambda}-d_{\lambda})\right\Vert
_{2}^{2}\leq\underline{u}(x_{\lambda},c_{\lambda})-\underline{u}(y_{\lambda
},d_{\lambda})+\overline{u}^{s_{0}}(x_{\lambda},c_{\lambda})-\overline
{u}^{s_{0}}(y_{\lambda},d_{\lambda})+4m(y_{\lambda}-x_{\lambda}).
\]
Consequently
\begin{equation}
\lambda\left\Vert (x_{\lambda}-y_{\lambda},c_{\lambda}-d_{\lambda})\right\Vert
_{2}^{2}\leq6m\left\Vert (x_{\lambda}-y_{\lambda},c_{\lambda}-d_{\lambda
})\right\Vert _{2}. \label{desxalfa-yalfac}%
\end{equation}
We can find a sequence $\lambda_{n}\rightarrow\infty$ such that $\left(
x_{\lambda_{n}},y_{\lambda_{n}},c_{\lambda_{n}},d_{\lambda_{n}}\right)
\rightarrow\left(  \widehat{x},\widehat{y},\widehat{c},\widehat{d}\right)  \in
A$. From (\ref{desxalfa-yalfac}), we get that
\begin{equation}
\left\Vert (x_{\lambda_{n}}-y_{\lambda_{n}},c_{\lambda_{n}}-d_{\lambda_{n}
})\right\Vert _{2}\leq6m/\lambda_{n} , \label{cota con 6m}%
\end{equation}
which gives $\widehat{x}=\widehat{y}$ and $\widehat{c}=\widehat{d}$.

Since $\Sigma^{\lambda}\left(  x,y,c,d\right)  =\underline{u}(x,c)-\overline
{u}^{s_{0}}(y,d)-\Phi^{\lambda}\left(  x,y,c,d\right)  $ reaches the maximum
in $\left(  x_{\lambda},y_{\lambda},c_{\lambda},d_{\lambda}\right)  \ $in the
interior of the set $A,$ the function%
\[
\psi(x,c)=\Phi^{\lambda}\left(  x,y_{\lambda},c,d_{\lambda}\right)
-\Phi^{\lambda}\left(  x_{\lambda},y_{\lambda},c_{\lambda},d_{\lambda}\right)
+\underline{u}\left(  x_{\lambda},c_{\lambda}\right)
\]
is a test for subsolution for $\underline{u}$ of the HJB equation at the point
$\left(  x_{\lambda},c_{\lambda}\right)  $. In addition, the function%
\[
\varphi^{s_{0}}(y,d)=-\Phi^{\lambda}\left(  x_{\lambda},y,c_{\lambda
},d\right)  +\Phi^{\lambda}\left(  x_{\lambda},y_{\lambda},c_{\lambda
},d_{\lambda}\right)  +\overline{u}^{s_{0}}\left(  y_{\lambda},d_{\lambda
}\right)
\]
is a test for supersolution for $\overline{u}^{s_{0}}$ at $\left(  y_{\lambda
},d_{\lambda}\right)  $ and so
\[
\partial_{c}\varphi^{s_{0}}(y_{\lambda},d_{\lambda})\leq-\frac{s_{0}}{c_{2}%
}\varphi(y_{\lambda},d_{\lambda})e^{-\frac{c}{c_{2}}}<0
\]
(because $y_{\lambda}>0)$. Hence, $\partial_{c}\psi(x_{\lambda},c_{\lambda
})=\partial_{c}\varphi^{s_{0}}(y_{\lambda},d_{\lambda})<0$, and we have
$\mathcal{L}^{c_{\lambda}}(\psi)(x_{\lambda},c_{\lambda})\geq0.$

\ Assume first that the functions $\underline{u}(x,c)$ and $\overline
{u}^{s_{0}}(y,d)$ are (2,1)-differentiable at $(x_{\lambda},c_{\lambda})\ $and
$(y_{\lambda},d_{\lambda})$ respectively. Since $\Sigma^{\lambda}$ defined in
(\ref{def-sigma-definitiva-por-hoy}) reaches a local maximum at $\left(
x_{\lambda},y_{\lambda},c_{\lambda},d_{\lambda}\right)  $ $\notin\partial A$,
we have that
\[
\partial_{x}\Sigma^{\lambda}\left(  x_{\lambda},y_{\lambda},c_{\lambda
},d_{\lambda}\right)  =\partial_{y}\Sigma^{\lambda}\left(  x_{\lambda
},y_{\lambda},c_{\lambda},d_{\lambda}\right)  =0
\]
and so%

\begin{equation}%
\begin{array}
[c]{lll}%
\partial_{x}\underline{u}(x_{\lambda},c_{\lambda}) & = & \partial_{x}
\Phi^{\lambda}(x_{\lambda},y_{\lambda},c_{\lambda},d_{\lambda})\\
& = & \lambda\left(  x_{\lambda}-y_{\lambda}\right)  +\frac{2m}{\left(
\lambda\left(  y_{\lambda}-x_{\lambda}\right)  +1\right)  ^{2}}\\
& = & -\partial_{y}\Phi^{\lambda}(x_{\lambda},y_{\lambda},c_{\lambda
},d_{\lambda})=\partial_{x}\overline{u}^{s_{0}}(y_{\lambda},d_{\lambda}).
\end{array}
\label{dos derivadas iguales}%
\end{equation}
Defining $A=\partial_{xx}\underline{u}(x_{\lambda},c_{\lambda})$ and
$B=\partial_{xx}\overline{u}^{s_{0}}(y_{\lambda},d_{\lambda})$, we obtain%

\[%
\begin{array}
[c]{l}%
\left(
\begin{array}
[c]{ll}%
\partial_{xx}\Sigma^{\lambda}\left(  x_{\lambda},y_{\lambda},c_{\lambda
},d_{\lambda}\right)  & \partial_{xy}\Sigma^{\lambda}\left(  x_{\lambda
},y_{\lambda},c_{\lambda},d_{\lambda}\right) \\
\partial_{xy}\Sigma^{\lambda}\left(  x_{\lambda},y_{\lambda},c_{\lambda
},d_{\lambda}\right)  & \partial_{yy}\Sigma^{\lambda}\left(  x_{\lambda
},y_{\lambda},c_{\lambda},d_{\lambda}\right)
\end{array}
\right) \\
=\left(
\begin{array}
[c]{ll}%
A-\partial_{xx}\Phi^{\lambda}\left(  x_{\lambda},y_{\lambda},c_{\lambda
},d_{\lambda}\right)  & -\partial_{xy}\Phi^{\lambda}\left(  x_{\lambda
},y_{\lambda},c_{\lambda},d_{\lambda}\right) \\
-\partial_{xy}\Phi^{\lambda}\left(  x_{\lambda},y_{\lambda},c_{\lambda
},d_{\lambda}\right)  & -B-\partial_{yy}\Phi^{\lambda}\left(  x_{\lambda
},y_{\lambda},c_{\lambda},d_{\lambda}\right)
\end{array}
\right)  \leq0.
\end{array}
\]
It is hence a negative semi-definite matrix, and%

\[%
\begin{pmatrix}
A & 0\\
0 & -B
\end{pmatrix}
\leq\partial_{xy}H\left(  \Phi^{\lambda}\right)  (x_{\lambda},y_{\lambda
},c_{\lambda},d_{\lambda}):=\left(
\begin{array}
[c]{ll}%
\partial_{xx}\Phi^{\lambda}\left(  x_{\lambda},y_{\lambda},c_{\lambda
},d_{\lambda}\right)  & \partial_{xy}\Phi^{\lambda}\left(  x_{\lambda
},y_{\lambda},c_{\lambda},d_{\lambda}\right) \\
\partial_{xy}\Phi^{\lambda}\left(  x_{\lambda},y_{\lambda},c_{\lambda
},d_{\lambda}\right)  & \partial_{yy}\Phi^{\lambda}\left(  x_{\lambda
},y_{\lambda},c_{\lambda},d_{\lambda}\right)
\end{array}
\right)  .
\]

In the case that $\underline{u}(x,c)$ and $\overline{u}^{s_{0}}(y,d)$ are not
(2,1)-differentiable at $\left(  x_{\lambda},c_{\lambda}\right)  \ $and
$(y_{\lambda},d_{\lambda})\ $respectively, we can resort to a more general
theorem to get a similar result. Using Theorem 3.2 of Crandall, Ishii and
Lions \cite{CrandallIshiLions}, it can be proved that for any $\delta>0$,
there exist real numbers $A_{\delta}$ and $B_{\delta}$ such that
\begin{equation}%
\begin{pmatrix}
A_{\delta} & 0\\
0 & -B_{\delta}%
\end{pmatrix}
\leq\partial_{xy}H\left(  \Phi^{\lambda}\right)  (x_{\lambda},y_{\lambda
},c_{\lambda},d_{\lambda})+\delta\left(  \partial_{xy}H\left(  \Phi^{\lambda
}\right)  (x_{\lambda},y_{\lambda},c_{\lambda},d_{\lambda})\right)  ^{2}
\label{MatrizDeltaruina}%
\end{equation}
and \bigskip%

\begin{equation}%
\begin{array}
[c]{c}%
\frac{\sigma^{2}}{2}A_{\delta}+(\mu-c_{\lambda})\partial_{x}\psi(x_{\lambda
},c_{\lambda})-q\psi(x_{\lambda},c_{\lambda})+c_{\lambda}\geq0,\\
\frac{\sigma^{2}}{2}B_{\delta}+(\mu-d_{\lambda})\partial_{x}\varphi^{s_{0}
}(y_{\lambda},d_{\lambda})-q\varphi^{s_{0}}(y_{\lambda},d_{\lambda
})+d_{\lambda}\leq0.
\end{array}
\label{Inecuaciones con A y B delta}%
\end{equation}
The expression (\ref{MatrizDeltaruina}) implies that $A_{\delta}-B_{\delta
}\leq0$ because%
\[
\partial_{xy}H\left(  \Phi^{\lambda}\right)  (x_{\lambda},y_{\lambda
},c_{\lambda},d_{\lambda})=\partial_{xx}\Phi^{\lambda}\left(  x_{\lambda
},y_{\lambda},c_{\lambda},d_{\lambda}\right)
\begin{pmatrix}
1 & -1\\
-1 & 1
\end{pmatrix}
\]
and
\[
\left(  \partial_{xy}H\left(  \Phi^{\lambda}\right)  (x_{\lambda},y_{\lambda
},c_{\lambda},d_{\lambda})\right)  ^{2}=2\left(  \partial_{xx}\Phi^{\lambda
}\left(  x_{\lambda},y_{\lambda},c_{\lambda},d_{\lambda}\right)  \right)
^{2}
\begin{pmatrix}
1 & -1\\
-1 & 1
\end{pmatrix}
.
\]
Therefore,%
\[%
\begin{array}
[c]{lll}%
A_{\delta}-B_{\delta} & = &
\begin{pmatrix}
1 & 1
\end{pmatrix}
\left(
\begin{array}
[c]{cc}%
A_{\delta} & 0\\
0 & -B_{\delta}%
\end{array}
\right)  \left(
\begin{array}
[c]{c}%
1\\
1
\end{array}
\right) \\
& \leq &
\begin{pmatrix}
1 & 1
\end{pmatrix}
\left(  \partial_{xy}H\left(  \Phi^{\lambda}\right)  (x_{\lambda},y_{\lambda
},c_{\lambda},d_{\lambda})+\delta\left(  \partial_{xy}H\left(  \Phi^{\lambda
}\right)  (x_{\lambda},y_{\lambda},c_{\lambda},d_{\lambda})\right)
^{2}\right)  \left(
\begin{array}
[c]{c}%
1\\
1
\end{array}
\right) \\
& = & 0.
\end{array}
\]
And so, since $\varphi^{s_{0}}\left(  y_{\lambda},d_{\lambda}\right)
=\overline{u}^{s_{0}}\left(  y_{\lambda},d_{\lambda}\right)  $, $\psi
(x_{\lambda},c_{\lambda})=\underline{u}(x_{\lambda},c_{\lambda})$ and
\[
\partial_{x}\varphi^{s_{0}}\left(  y_{\lambda},d_{\lambda}\right)
=-\partial_{y}\Phi^{\lambda}\left(  x_{\lambda},y_{\lambda},c_{\lambda
},d_{\lambda}\right)  =\partial_{x}\Phi^{\lambda}\left(  x_{\lambda
},y_{\lambda},c_{\lambda},d_{\lambda}\right)  =\partial_{x}\psi(x_{\lambda
},c_{\lambda}),
\]
we obtain%

\begin{equation}%
\begin{array}
[c]{lll}%
\underline{u}(x_{\lambda},c_{\lambda})-\overline{u}^{s_{0}}\left(  y_{\lambda
},d_{\lambda}\right)  & = & \psi(x_{\lambda},c_{\lambda})-\varphi^{s_{0}
}\left(  y_{\lambda},d_{\lambda}\right) \\
& \leq & \frac{\sigma^{2}}{2q}(A_{\delta}-B_{\delta})\\
&  & +\left(  \frac{c_{\lambda}}{q}-\frac{d_{\lambda}}{q}\right)
(1-\partial_{x}\Phi^{\lambda}\left(  x_{\lambda},y_{\lambda},c_{\lambda
},d_{\lambda}\right)  )\\
& \leq & \left(  \frac{c_{\lambda}}{q}-\frac{d_{\lambda}}{q}\right)
(1-\lambda\left(  x_{\lambda}-y_{\lambda}\right)  -\frac{2m}{\left(
\lambda\left(  y_{\lambda}-x_{\lambda}\right)  +1\right)  ^{2}}).
\end{array}
\label{desdifsup1}%
\end{equation}
Hence, from (\ref{cota con 6m}) and (\ref{liminf-mlambda}), we get%

\begin{align*}
0  &  <M\leq\liminf\limits_{\lambda\rightarrow\infty}M_{\lambda}\leq
\lim\limits_{_{n\rightarrow\infty}}M_{\lambda_{n}}=\lim\limits_{_{n\rightarrow
\infty}}\Sigma^{\lambda_{n}}(x_{\lambda_{n}},y_{\lambda_{n}},c_{\lambda_{n}
},d_{\lambda_{n}})=\underline{u}(\widehat{x},\widehat{c})-\overline{u}^{s_{0}
}(\widehat{x},\widehat{c})\\
&  \leq\lim_{n\longrightarrow\infty}\left(  \frac{c_{\lambda_{n}}}{q}
-\frac{d_{\lambda_{n}}}{q}\right)  (1-\lambda_{n}\left(  x_{\lambda_{n}
}-y_{\lambda_{n}}\right)  -\frac{2m}{\left(  \lambda_{n}\left(  y_{\lambda
_{n}}-x_{\lambda_{n}}\right)  +1\right)  ^{2}})\\
&  \leq.\lim_{n\longrightarrow\infty}\left\vert \frac{c_{\lambda_{n}}}
{q}-\frac{d_{\lambda_{n}}}{q}\right\vert (1+\lambda_{n}\left\Vert
(x_{\lambda_{n}}-y_{\lambda_{n}},c_{\lambda_{n}}-d_{\lambda_{n}})\right\Vert
_{2}+\frac{2m}{\left(  \lambda_{n}\left(  y_{\lambda_{n}}-x_{\lambda_{n}
}\right)  +1\right)  ^{2}})\\
&  \leq\lim_{n\longrightarrow\infty}\left\vert \frac{c_{\lambda_{n}}}{q}
-\frac{d_{\lambda_{n}}}{q}\right\vert (1+8m)=0.
\end{align*}
This is a contradiction and so we get the result. $\blacksquare$\newline

\end{document}